\documentclass[reqno]{amsart}
\usepackage[process=all,crop=pdfcrop]{pstool}
\usepackage{float}
\usepackage{pdfpages}
\numberwithin{equation}{section} % Equation numbering control.
\usepackage{comment}
\usepackage{geometry}
\usepackage{amsmath}
\usepackage{amssymb}
\usepackage{graphicx}
\usepackage{soul}
\usepackage{enumitem}
\usepackage{mathptmx}  
\usepackage[utf8]{inputenc}
\usepackage{mathtools}
\usepackage{lipsum}
\usepackage{color}
\usepackage{subfigure}
\usepackage{caption}
\usepackage{wrapfig}
\usepackage{xcolor}
\usepackage[rightcaption]{sidecap}
\usepackage[square,comma,numbers,sort&compress]{natbib}
\usepackage[colorlinks=true]{hyperref}
\hypersetup{urlcolor=blue, citecolor=red}

\DeclareGraphicsExtensions{.eps,.pdf,.png,.jpg}

\newcommand{\Revtwo}[1]{{\color{red} #1}}

\begin{document}

\title{Time-Dependent Duhamel Renormalization method with Multiple Conservation and Dissipation Laws }


%
\author{Sathyanarayanan Chandramouli}
\address[S. Chandramouli]{Department of Mathematics\\
                Florida State University\\
       Tallahassee, FL 32306, USA}
\email{schandra@math.fsu.edu}

\author{Aseel Farhat}
\address[A. Farhat]{Department of Mathematics\\
                Florida State University\\
       Tallahassee, FL 32306, USA}
\email{afarhat@fsu.edu}

\author{Ziad Musslimani}
\address[Z. Musslimani]{Department of Mathematics\\
                Florida State University\\
       Tallahassee, FL 32306, USA}
\email{musliman@math.fsu.edu}

%%%%%%%%%%%%%%%%%%%%%%%%%%%%%%%%%%%%%%%%%%%%%%%%%%%%%%%%%%%

\begin{abstract}

The time dependent spectral renormalization (TDSR) method was introduced by Cole and Musslimani as a novel way to numerically solve initial boundary value problems. An important and novel aspect of the TDSR scheme is its ability to incorporate physics in the form of conservation laws or dissipation rate equations. However, the method was limited to include a single conserved or dissipative quantity.

The present work significantly extends the computational features of the method with the (i) incorporation of multiple conservation laws and/or dissipation rate equations, (ii) ability to enforce versatile boundary conditions, and (iii) higher order time integration strategy. The TDSR method is applied on several prototypical evolution equations of physical significance. Examples include the Korteweg-de Vries (KdV), multi-dimensional nonlinear Schr\"odinger (NLS) and the Allen-Cahn equations. 

\end{abstract}

\maketitle

{\bf Keywords:} Renormalization method, initial boundary value problems, partial differential equations, Duhamel's principle, nonlinear waves, soliton equations, Hamiltonian and dissipative systems.

\bigskip 
%%%%%%%%%%%%%%%%%%%%%%%%%%%%%%%%%%%%%
\section{Introduction}
%%%%%%%%%%%%%%%%%%%%%%%%%%%%%%%%%%%%%%
Numerical simulation of initial boundary value problems is of utmost importance in several engineering and scientific disciplines. Over the last few decades, several time-stepping methods have been developed and proposed to achieve this goal. Among them are the class of implicit/explicit Runge-Kutta methods \cite{quarteroni2010numerical}, exponential time-differencing \cite{taflove1998advances,cox2002exponential,kassam2005fourth,yang2020fidelity} and the split-step operator splitting  \cite{wiedeman1986} to name a few. For evolution equations that arise in physical applications, it is highly desirable to devise time-stepping schemes that reflect the underlying physics at hand.
Such structure preserving numerical schemes are of paramount importance for long-time integration, where either it is necessary to ensure numerical stability (e.g., if the numerics could conserve the $L^2$ norm of the solution for the KdV/NLS), or preserve other features ({such as capturing the correct shock speed in the context of systems of hyperbolic partial differential equations}).% designing conservative numerical methods for the discretization of system of hyperbolic equations is necessary to ensure the correct shock speed }). 
 To date, there are various ways to input some physics into the numerical time-integration. For example, the geometric/symplectic integrators that preserve the Hamiltonian and symplectic structure \cite{hairer2006geometric}, the operator splitting method that was used for the NLS equation to preserve the power and the non-linear dispersion relation \cite{wiedeman1986}, the multi-symplectic  schemes designed for the generalized Schr{\"o}dinger equations \cite{ISLAS2003403,islas2001geometric,ISLAS2005290}, the Taha-Ablowitz \cite{taha1984analytical}, the Ablowitz-Ladik \cite{ablowitz1976nonlinear}, and the Ablowitz-Musslimani \cite{ablowitz2014integrable,ablowitz2020discrete} schemes that preserve the integrable structure of the KdV, NLS and class of nonlocal NLS equations, respectively. Other relevant works  correspond to {the conservative finite volume Godunov schemes
%used to discretize the compressible gas dynamics (Euler) equations 
\cite{godunov1959difference,hirsch1990numerical, Leveque}}, finite difference schemes that preserve the energy or  dissipation property of the underlying model equation (see, e.g.,\cite{FURIHATA1999181,sanz1982explicit}), {as well as} a finite volume scheme that {conserve} mass and momentum of the KdV equation (see, e.g., \cite{cui2007numerical}).
Recently, Cole and Musslimani \cite{cole2017time} proposed an alternative method to simulate dynamical systems that enables the inclusion of physics ``on-demand". The core idea is to make use of the Duhamel's principle to recast the underlying evolution equation as a space-time integral equation. The resulting system is then solved iteratively using a novel time-dependent renormalization process that controls both the numerical convergence properties of the scheme while at the same time preserving a single physical law. 

In this paper, we extend the time-dependent spectral renormalization method to allow for multiple conservation laws or dissipation rates to be simultaneously incorporated in the simulation. This is achieved by introducing as many time-dependent renormalization factors as the number of conservation/dissipation laws being enforced. The solution sought is then written as a linear superposition of finite number of space-time dependent auxiliary wave functions with the time-dependent renormalization factors envisioned to play the role of ``expansion coefficients". When inserted into the corresponding Duhamel's formula, a finite set of ``sub-Duhamel" integral equations are obtained governing the space-time dynamics of each individual auxiliary function. These integral equations are then numerically solved using a novel space-time fixed point iteration. The importance of such a dynamic renormalization process are twofold: (i) it provides convergence when needed and (ii) enables the inclusion of conservation/dissipation laws.  
The Duhamel integrals are numerically computed using a Cauchy-Filon-Simpson quadrature formula.
The TDSR method is implemented on several prototypical evolution equations of physical significance. This includes the KdV,  Allen-Cahn, multi-dimensional and the $PT$ symmetric integrable nonlocal NLS equations. 

%%%%%%%%%%%%%%%%%%%%%%%%%%%%%%%%%%%%%%%%%%%%%%%%%%%%%%%%%%%
The paper is organized as follows. In Sec.~\ref{Sec-TDSR}  we put forward a general framework for the TDSR scheme in arbitrary space dimensions and show how to incorporate multiple conservation laws or dissipation rate equations into the algorithm. The Cauchy-Filon Simpson time integration is derived in Sec.~\ref{time-int-periodic} for evolution equations subject to either periodic or rapidly decaying boundary conditions. The TDSR scheme is applied on the KdV and NLS equations, with single and multiple conservation laws.
In Sec.~\ref{time-integration-non-periodic} the Cauchy-Filon trapezoidal time integration is developed for evolution equations subject to non-periodic and non decaying boundary conditions. In this regard, the Allen-Cahn equation is used as a test bed to assess the performance of the algorithm. We conclude in Sec.~\ref{sec-conclusions} with comments on future directions.

\bigskip 
\bigskip
%%%%%%%%%%%%%%%%%%%%%%%%%%%%%%%%%%%%%%%%%%%%%%%%%%%%%%%%%%
\section{TDSR and Duhamel principle}
\label{Sec-TDSR}
%%%%%%%%%%%%%%%%%%%%%%%%%%%%%%%%%%%%%%%%%%%%%%%%%%%%%%%%%%

In this section, we formulate the TDSR method using Duhamel's principle in conjunction with multiple conservation laws. Consider the evolution equation for the real (or complex) valued function $u(\textbf{x},t)$:
\begin{equation}
\label{pde}
u_t  = \mathcal{L}(u) + \mathcal{N}(u), \;\;\;\; u({\bf x}, 0) = u_0({\bf x}),
\end{equation}
where $\mathcal{L}$ is a linear, constant coefficient differential operator
and $\mathcal{N}(u)$ is a nonlinear {operator}.  The initial-boundary value problem (\ref{pde}) is posed on a spatial domain $\Omega$ that is either bounded or unbounded. Furthermore, Eq.~(\ref{pde}) is supplemented with either periodic, rapidly decaying, or other types of boundary conditions. 
As mentioned above, we are interested in evolution equations that are either (i) conservative, in which case, there exists $N$ conserved quantities given by
 \begin{equation}
 \label{n-cons-laws}
 \mathcal{Q}_m(u)\equiv
 \int_{\Omega} Q_{m}[u(\textbf{x},t)] d{ \textbf{x}} = \int_{\Omega} Q_{m}[u_0(\textbf{x})] d{\bf x}\equiv C_{m}\;,\;\;\;\; m=1,2,3,\cdots N,
 \end{equation}
or (ii) dissipative, so that there are $N$ densities $\rho_m$ and fluxes $F_m$ that obey the rate equations
\begin{equation}
 \label{n-diss-rates}
        \frac{d}{dt} {\int_{\Omega} \rho_m[u(\textbf{x},t)]d{\bf x}} = - {\int_{\Omega} F_m[u(\textbf{x},t)] d{\bf x}} \;,\;\;\; m=1,2,3,\cdots , N.
    \end{equation}
    %%%%%%%%%%%%%%%%%%%%%%%%%%%%%%%%%%%%%%%%%%%%%%%%%%%%%%%
Equation~(\ref{pde}) is rewritten in an integral form using the Duhamel's principle:
%%%%%%%%%%%%%%%%%%%%%%%%%%%%%%%%%%%%%%%%%%%%%%%%%%%%%%%%%%
\begin{align}
\label{Duhamel}
  \phantom{u(x,t)}
  &\begin{aligned}
    \mathllap{u(\textbf{x},t)} & = e^{t\mathcal{L}}[u_0(\textbf{x})]+ \int_{0}^{t} e^{(t-\tau) \mathcal{L}}\mathcal{N}[u({\bf x},\tau) ]d{\tau}.
  \end{aligned}
\end{align}
%%%%%%%%%%%%%%%%%%%%%%%%%%%%%%%%%%%%%%%%%%%%%%%%%%%%%%%%%%%%%%%%%%%%%%%%%%%%
Recall that for any periodic or $L^2$ function $w({\bf x})$, the semi-group $e^{t\mathcal{L}}$ admits the spectral representation:
\begin{equation}
\label{Linear-operator-semigroup}
e^{t\mathcal{L}}[w(\textbf{x})] = \mathcal{F}^{-1} [\exp (t\hat{\mathcal{L}}) \mathcal{F} [w(\textbf{x})]],
\end{equation}
where $\hat{\mathcal{L}}$ is the Fourier symbol associated with $\mathcal{L}$ and $\mathcal{F}$, $\mathcal{F}^{-1}$ denote the $d$-dimensional forward and inverse Fourier transforms, defined by
%%%%%%%%%%%%%%%%
\begin{eqnarray}
     \hat{w}(\textbf{k})  \equiv \mathcal{F}[w(\textbf{x})] =
    (2\pi)^{-d/2}\int_{\mathbb{R}^d} w(\textbf{x})e^{-i\textbf{k} \cdot \textbf{x}}d{\textbf{x}},  \label{Forward-trans}\\
    \mathcal{F}^{-1}[\hat{w}(\textbf{k})]=(2\pi)^{-d/2}\int_{\mathbb{R}^d}\hat{w}(\textbf{k})e^{i\textbf{k}\cdot\textbf{x}}d{\textbf{k}}.
\end{eqnarray}
%%%%%%%%%%%%%%%%%
For periodic functions defined on a bounded spatial domain the forward Fourier integral \eqref{Forward-trans} represents the coefficients of its Fourier series.

We now outline in details how to incorporate multiple conservation laws or dissipation rate equations in the TDSR scheme. To this end, we seek a solution to Eq.~(\ref{Duhamel}) in the form
%%%%%%%%%%%%%%%%%%%%%%%%%%%%%%%%%
     \begin{equation}
     \label{decomposition}
         u(\textbf{x},t) = \sum_{j=1}^{N} R_{j}(t) v_{j}(\textbf{x},t),
     \end{equation}
%%%%%%%%%%%%%%%%%%%%%%%%%%%%%%%%%     
where $R_j(t)$ are time-dependent renormalization factors to be determined from knowledge of the conservation or dissipation laws and $v_j(\textbf{x},t)$ are space-time dependent auxiliary functions that satisfy the same boundary conditions as the solution $ u(\textbf{x},t)$. Our extensive numerical experiments seem to indicate that in order for the TDSR iteration to converge, the initial condition $u_0(\textbf{x})$ needs to be re-written as a sum 
of $N$ ({identically non-zero}) functions $f_j(\textbf{x})$ (here referred to as pseudo initial conditions). Namely, we write
 %%%%%%%%%%%%%%%%
   \begin{equation}
   \label{fj-condition}
 u_0(\textbf{x})   =  \sum_{j=1}^N f_j(\textbf{x}),
   \end{equation}
%%%%%%%%%%%%%%%%%%%%    
where each $f_j$ is chosen to be compatible with the underlying boundary conditions. The specific choice of the functions $\{f_j({\bf x})\}$, $j=1,2,\cdots,N$, is discussed in Sec.~(\ref{Multiple-conservation-laws-section}) when solving the 
KdV equation.
Substituting Eqns.~(\ref{decomposition}) and (\ref{fj-condition}) into (\ref{Duhamel}), we obtain an equation for the auxiliary functions $v_{j}(\textbf{x},t), j=1,2,\cdots N$:
%%%%%%%%%%%%%%%%%%%%%%%%%%%%%%%%%%%%%%%%%%%
\begin{equation}
\label{Duhamel-2}
 \sum_{j=1}^{N} R_{j}(t) v_{j}(\textbf{x},t)  =   \sum_{j=1}^N e^{t\mathcal{L}}[f_j(\textbf{x}) ]
 + 
 \int_{0}^{t} e^{(t-\tau) \mathcal{L}}\mathcal{N} \left[  \sum_{j=1}^{N} R_{j}(\tau) v_{j}(\textbf{x},\tau ) \right]d{\tau}.
\end{equation}
%%%%%%%%%%%%%%%%%%%%%%%%%%%%%%%%%%%%%%%%%%%
Scrutinizing Eq.~(\ref{Duhamel-2}) reveals the existence of $N-1$ degrees of freedom for the variables $v_1, v_2,\cdots, v_{N-1}.$ Next, we outline how to eliminate each degree of freedom and derive a self-consistent set of equations that forms the basis for the TDSR scheme. To begin with, we choose $v_{1}(\textbf{x},t)$ such that
%%%%%%%%%%%%%%%%%%%%%%%%%%%
     \begin{equation}
     \label{eqn-phi-1}
 v_{1}(\textbf{x},t) = \mathcal{M}_1[R_1,v_1] \equiv
 \frac{1}{R_1(t)}\left\{ e^{t\mathcal{L}}[f_1(\textbf{x}) ]+\int_{0}^{t} e^{(t-\tau)\mathcal{L}}\mathcal{N}[ R_1(\tau ) v_{1}(\textbf{x},\tau )]d{\tau} \right\}.
     \end{equation}
%%%%%%%%%%%%%%%%%%%%%%%%%%%
The rationale behind this choice is rooted in the fact that 
Eq.~(\ref{eqn-phi-1}) must reduce back to the case when only one conservation or dissipation law is under consideration with 
$f_1(\textbf{x}) \equiv u_0(\textbf{x}).$  With this at hand, we next require $v_2(\textbf{x},t)$ to satisfy
%%%%%%%%%%%%%%%%%%%%%%%%%%%
%%%%%%%%%%%%%%%%%%%%%%%%%%%
{
\begin{eqnarray}
\label{eqn-phi-2}
v_2(\textbf{x},t)&=&\mathcal{M}_2[R_1,R_2,v_1,v_2] \nonumber\\&\equiv&  \frac{1}{R_2(t)} \Big(e^{t\mathcal{L}}[f_2(\textbf{x}) ] 
% \nonumber \\
 +  \int_{0}^{t} e^{(t-\tau)\mathcal{L}} \mathcal{N} \left[ R_1 (\tau ) v_{1}(\textbf{x},\tau )+ R_2(\tau )v_{2}(\textbf{x},\tau )\right] d{\tau}\Big)
 \nonumber \\
 &-&  \frac{1}{R_2(t)} \int_{0}^{t} e^{(t-\tau)\mathcal{L}}  \mathcal{N} \left[R_1 (\tau ) v_{1}(\textbf{x},\tau ) \right]  d{\tau},\\\nonumber 
\end{eqnarray}
and for general $j=3,4,\cdots, N$, 
}
{
     \begin{eqnarray}
     \label{eqn-phi-j}
 v_{j}(\textbf{x},t) &=&   \mathcal{M}_j[R_1, R_2,\cdots, R_j;v_1,v_2,\cdots, v_j]\\\nonumber
  &\equiv& \frac{1}{R_j(t)} e^{t\mathcal{L}}[f_j(\textbf{x}) ] 
% \nonumber \\
 +  \frac{1}{R_j(t)} \int_{0}^{t} e^{(t-\tau)\mathcal{L}} \mathcal{N} \left[ \sum_{\ell =1}^{j} R_\ell (\tau ) v_{\ell}(\textbf{x},\tau ) \right] d{\tau}
  \\\nonumber
 &-&  \frac{1}{R_j(t)} \int_{0}^{t} e^{(t-\tau)\mathcal{L}}  \mathcal{N} \left[ \sum_{\ell =1}^{j-1} R_\ell (\tau ) v_{\ell}(\textbf{x},\tau ) \right]  d{\tau},
     \end{eqnarray}}
%%%%%%%%%%%%%%%%%%%%%%%%%%%
%where $ \Revtwo{\mathcal{M}_j[R,v] \equiv \mathcal{M}_j[R_1, R_2,\cdots, R_j;v_1,v_2,\cdots, v_j].}$
Note that Eqns.~(\ref{eqn-phi-1})-(\ref{eqn-phi-j}) are self-consistent with the Duhamel's formula (\ref{Duhamel-2}). Indeed, multiplying (\ref{eqn-phi-2}) by $R_{2}$ and (\ref{eqn-phi-j}) by $R_j$ and summing over all $j= 2, 3,\cdots, N$, we obtain
%%%%%%%%%%%%%%%%%%%%%%%%%
%%%%%%%%%%%%%%%%%%%%%%%%%%%%%%%%%%%
\begin{eqnarray}
     \label{eqn-phi-j3}
\sum_{j=2}^{N} R_{j}(t) v_{j}(\textbf{x},t) &=&  e^{t\mathcal{L}}\left[\sum_{j=2}^{N}f_j(\textbf{x}) \right] 
 \nonumber \\
 &+& \int_{0}^{t} e^{(t-\tau)\mathcal{L}}  \sum_{j=2}^{N} \left\{ \mathcal{N} \left[ \sum_{\ell =1}^{j} R_\ell (\tau ) v_{\ell}(\textbf{x},\tau ) \right] \right\} d{\tau}
 \nonumber \\
 &-&   \int_{0}^{t} e^{(t-\tau)\mathcal{L}}  \sum_{j=2}^{N} \left\{ \mathcal{N} \left[ \sum_{\ell =1}^{j-1} R_\ell (\tau ) v_{\ell}(\textbf{x},\tau ) \right] \right\} d{\tau}.
\end{eqnarray}
%%%%%%%%%%%%%%%%%%%%%%%%%%%%%%%%%%%%
The last two terms on the right hand side of Eq.~(\ref{eqn-phi-j3}) satisfy
%%%%%%%%%%%%%%%%%%%%%%%%%%%%%%%%%%%%
\begin{eqnarray}
     \label{eqn-telescopic-sum}
 & \int_{0}^{t} d\tau e^{(t-\tau)\mathcal{L}}
 \sum_{j=2}^{N} \left\{ \mathcal{N}\left[ \sum_{\ell =1}^{j} R_\ell (\tau ) v_{\ell}(\textbf{x},\tau ) \right]
-   \mathcal{N} \left[ \sum_{\ell =1}^{j-1} R_\ell (\tau ) v_{\ell}(\textbf{x},\tau ) \right]  \right\}
  \nonumber \\
&\qquad \qquad =  
 \int_{0}^{t} d\tau e^{(t-\tau)\mathcal{L}}
\mathcal{N} \left[ \sum_{\ell =1}^{N} R_\ell (\tau ) v_{\ell}(\textbf{x},\tau ) \right] 
- \mathcal{N}[ R_1(\tau ) v_{1}(\textbf{x},\tau )].
\end{eqnarray}
%%%%%%%%%%%%%%%%%%%%%%%%%%%%%%%%%%%%%
The conclusion is complete once we multiply Eq.~(\ref{eqn-phi-1}) by $R_{1}(t);$ add the result to 
Eq.~(\ref{eqn-telescopic-sum}) and use the condition (\ref{fj-condition}). In summary, Eqns.~(\ref{eqn-phi-1})-(\ref{eqn-phi-j}) give an implicit integral representation for the auxiliary functions $ v_j.$ To close the system,
all we need is to compute the renormalization factors $R_j(t)$. Substituting Eq.~(\ref{decomposition}) into 
Eqns.~(\ref{n-cons-laws}) and (\ref{n-diss-rates}) gives:\\\\
%%%%%%%%%%%%%%%%%%%%%%%%%%%%%%%%%%%%%%%%%%%%%%%%%%%%%%%
{\bf Conservative case:}
     \begin{equation}
     \label{cons-laws-R}
     \mathcal{Q}_m \left( \sum_{\ell=1}^{N} R_{\ell}(t) v_{\ell}({\bf x},t) \right)   \equiv
 \int_{\Omega} Q_{m}\left( \sum_{\ell=1}^{N} R_{\ell}(t) v_{\ell}({\bf x},t) \right)d {\bf x} = C_{m},
     \end{equation}
%%%%%%%%%%%%
{\bf Dissipative case:}
%%%%%%%%%%%%%%%%%%%%%%%%%%%%%%%%%%%%%%%%%%%%%%%%%%%%%%%%%%%%%
\begin{equation}
 \label{n-diss-rates-R}
   \frac{d}{dt} {\int_{\Omega} \rho_m\left[ \sum_{\ell =1}^{N} R_{\ell}(t) v_{\ell}({\bf x},t)  \right] d{\bf x}} 
        = - {\int_{\Omega} F_m \left[ \sum_{\ell =1}^{N} R_{\ell}(t) v_{\ell}({\bf x},t)  \right] d{\bf x}},
            \end{equation}
where $ m=1,2,\cdots , N$. System (\ref{cons-laws-R}) defines $N$ algebraic equations for the time-dependent renormalization factors whereas (\ref{n-diss-rates-R}) a set of coupled nonlinear ordinary differential equations governing the evolution of $R_j(t)$. With this at hand, the TDSR iterative process is summarized below:
%%%%%%%%%%%%%%%%%%%%%%%%%%%
     \begin{align}
     \label{eqn-phi-iteration}
 &v^{(n+1)}_{1} = \mathcal{M}_1[R_1^{(n)}(t),v_1^{(n)}],\\
     %\end{equation}
%%%%%%%%%%%%%%%%%%%%%%%%%%%
%%%%%%%%%%%%%%%%%%%%%%%%%%%
\label{eqn-phi-iteration-2}
 &v^{(n+1)}_{2} = \mathcal{M}_2[R_1^{(n)},R_2^{(n)},v_1^{(n)},v_2^{(n)}],\\
    % \begin{equation}
     \label{eqn-phi-iteration-j}
 &v^{(n+1)}_{j} = \mathcal{M}_j[R_1^{(n)},R_2^{(n)},\cdots , R_j^{(n)},v_1^{(n)},v_2^{(n)},\cdots , v_j^{(n)}]\;,\; j=3,\cdots, N,
     \end{align}
%%%%%%%%%%%%%%%%%%%%%%%%%%%
%%%%%%%%%%%%%%%%%%%%%%%%%%%
 with $R_j^{(n)}\;,j=1,2,\cdots , N$ given by  (for conservative cases) 
\begin{align}
 \mathcal{Q}_m \left( \sum_{\ell=1}^{N} R^{(n)}_{\ell}(t) v^{(n)}_{\ell}({\bf x},t) \right)     \equiv C_m\;,
\label{eqn-Q-iteration}
\end{align}
and
\begin{align}
   \frac{d}{dt} {\int_{\Omega} \rho_m\left[ \sum_{\ell =1}^{N} R^{(n)}_{\ell}(t) v^{(n)}_{\ell}({\bf x},t)  \right] d{\bf x}} 
        = - {\int_{\Omega} F_m \left[ \sum_{\ell =1}^{N} R^{(n)}_{\ell}(t) v^{(n)}_{\ell}({\bf x},t)  \right] d{\bf x}}, \label{eqn-Y-iteration}
\end{align}
for the dissipative cases where $ m=1,2,\cdots, N.$ As a reminder, the functionals $\mathcal{M}_1$ and $\mathcal{M}_j,  j=2,3,\cdots , N$ are respectively defined by Eqns.~(\ref{eqn-phi-1}) \Revtwo{-} (\ref{eqn-phi-j}). A workflow for the TDSR algorithm is given below, clarifying the structure of the iterative process:
%%%%%%%%%%%%%%%%%%%%%%%%%%%%%%%%%%%%%%%%%%%%%%%%%%%%%%%%%%%%%%%%%%%%%%%%%
\begin{enumerate}
\item {\bf Choose the pseudo initial conditions}: The set of pseudo initial conditions $f_j({\bf x})$, $j=1,2,\cdots\;N$, are chosen in such a way that Eq.~\eqref{fj-condition} is satisfied (see Sec.\ref{Multiple-conservation-laws-section} for further details). Note that they are used in the Picard iterations defined by Eqs.~(\ref{eqn-phi-1})-(\ref{eqn-phi-j}). 
    %%%%%%%%%%%%%%%%%%%%%%%%%%%%%%%%%%%%%%%%%%%%%%%%%%%%%%%%%%%%%%%%%%%%%%%%%
 \item {\bf Select initial guesses $v_j^{(1)}(\textbf{x},t)$ for $j=1,2,\cdots, N$}:
We seed Eqs.~(\ref{eqn-phi-1})-(\ref{eqn-phi-j})  with these initial guesses for the space-time dependent auxiliary functions. 
%%%%%%%%%%%%%%%%%%%%%%%%%%%%%%%%%%%%%%%%%%%%%%%%%%%%%%%%%%%%%%%%%%%%%%%%%%%
\item {\bf Compute the initial iterate of the set of renormalization factors}:
The set of auxiliary functions are used to evaluate
 the time-dependent renormalization factors ${R_j^{(1)}(t)}$ , $j=1,2,\cdots,N$ via the system of equations (\ref{eqn-Q-iteration}) or
(\ref{eqn-Y-iteration}) depending on whether the underlying evolution equation is conservative or dissipative.
%%%%%%%%%%%%%%%%%%%%%%%%%%%%%%%%%%%%%%%%%%%%%%%%%%%%%%%%%%%%%%%%%%%%%%%%%
\item {\bf Compute the Duhamel integrals defined in Eqs.~\eqref{eqn-phi-iteration}-\eqref{eqn-phi-iteration-j}}: The Duhamel integrals are computed using $v_j^{(1)}({\bf x},t)$ and $R_j^{(1)}(t)$, for $j=1,2,\cdots,N$. 
%%%%%%%%%%%%%%%%%%%%%%%%%%%%%%%%%%%%%%%%%%%%%%%%%%%%%%%%%%%%%%%%%%%%%%%%%%%
\item {\bf Update the Duhamel iteration}:
The Duhamel integrals computed in the previous step are now used to compute the second iterate of $v_j^{(2)}({\bf x},t)$ using Eqs.~\eqref{eqn-phi-iteration}-\eqref{eqn-phi-iteration-j}.
%%%%%%%%%%%%%%%%%%%%%%%%%%%%%%%%%%%%%%%%%%%%%%%%%%%%%%%%%%%%%%%%%%%%%%%%%%%%
\item {\bf Update the renormalization factors}: The updated $\{v_j^{(2)}({\bf x},t)\}$'s are now used to correct the the set of renormalization factors $\{R_j^{(2)}(t)\}$ , $j=1,2,\cdots,N$; from the system of equations given by \Revtwo{\eqref{eqn-Q-iteration}} (for the conservative case), or \Revtwo{\eqref{eqn-Y-iteration}} (for the dissipative case). 
\item {\bf Iterative update}: Repeat steps (5) and (6) till convergence is achieved.
%%%%%%%%%%%%%%%%%%%%%%%%%%%%%%%%%%%%%%%%%%%%%%%%%%%%%%%%%%%%%%%%%%%%%%%%%%%
% \item {\bf The iterative process}  The TDSR iterations are comprised of the following steps:
 %   \begin{enumerate}
  %      \item  Compute $v_j^{(n+1)}({\bf x},t)$, for $j=1,2,\cdots,N$ using Eqs.(\ref{eqn-phi-iteration})-(\ref{eqn-phi-iteration-j}). The Duhamel integrals are evaluated from the set of renormalization factors $\{R_j^{(n)}(t)\}$ and set of space-time dependent auxiliary functions $\{v_j^{(n)}(x,t)\}$ output at the $(n)^{th}$ TDSR iteration.
   %     \item FIX THE SENTENCE: The updated $\{v_j^{(n+1)}({\bf x},t)\}$'s are then used to correct the the set of time-dependent renormalization factors $\{R_j^{(n+1)}(t)\}$ , $j=1,2,\cdots,N$; from the system of equation \Revtwo{(\ref{eqn-Q-iteration})} for the conservative case, or the system of equations in \Revtwo{(\ref{eqn-Y-iteration})} for the dissipative case. This in turn ensures that the conservation laws/dissipation rate equations are enforced at {every}  iteration. The process is repeated till convergence is achieved.
    %\end{enumerate}
\end{enumerate}
%%%%%%%%%%%%%%%%%%%%%%%%%%%%%%%%%%%%%%%%%%%%%%%%%%%%%%%%%%%%%

\bigskip
\bigskip 
%%%%%%%%%%%%%%%%%%%%%%%%%%%%%%%%%%%%%%%%%%%%%%%%%%%%%%
\section{Time integration with various boundary conditions}
\label{time_int}
%%%%%%%%%%%%%%%%%%%%%%%%%%%%%%%%%%%%%%%%%%%%%%%%%%%
\subsection{Periodic and decaying boundary conditions}\label{time-int-periodic}
 %%%%%%%%%%%%%%%%%%%%%%%%%%%%%%%%%%%%%%%%%%%%%%%%
 %%%%%%%%%%%%%%%%%%%%%%%%%%%%%%%%%%%%%%%%%%%%%%%%
 In this section, we detail the numerical approach used to approximate the Duhamel integral
\begin{equation}
  \label{time-integral}
 I(\textbf{x},t) \equiv \int_{0}^{t} e^{(t-\tau) \mathcal{L}} G(\textbf{x},\tau )d \tau,
  \end{equation}
with $G(\textbf{x},\tau ) \equiv \mathcal{N}[u(\textbf{x},\tau )].$ When subject to periodic or rapidly decaying boundary conditions, the action of the semi-group ${\exp}(t\mathcal{L})$ on $G$ follows from its spectral representation $\mathcal{F}[\exp(t \mathcal{L}) G]=\exp [t \hat{\mathcal{L}}]\mathcal{F}[G]$ where $\hat{\mathcal{L}}$ is the Fourier symbol associated with the constant coefficients linear operator $\mathcal{L}.$ 
Our approach in approximating the integral Eq.~\eqref{time-integral} is based on the Filon-Simpson quadrature method  \cite{iserles2004numerical,olver2008numerical,ascher2011first}. To this end, we consider an $N_T$ equally spaced mesh points residing inside the time interval $[0,T]$ with
$t_i=i\Delta t, i=0,1,2,\cdots, N_T$, labeling all grid points. It can be shown that $\hat{I}(\textbf{k},t_{i})$ {satisfies} the exact recurrence relation
	%%%%%%%%%%%%%%%%%
	\begin{equation}
	\label{I-simpson-int}
	\hat{I}(\textbf{k},t_{i+1})= e^{2\Delta t\hat{\mathcal{L}}(\textbf{k})} \left[\hat{I}(\textbf{k},t_{i-1})+ \int_{t_{i-1}}^{t_{i+1}}e^{(t_{i-1}-\tau)\hat{\mathcal{L}}(\textbf{k})}\hat{G}(\textbf{k},\tau)d{\tau} \right].
	\end{equation} 
	%%%%%%%%%%%%%%%%%
As a reminder, a hat over a quantity represents its Fourier transform (see definition (\ref{Forward-trans})) or its Fourier series coefficients. 
Next, we approximate the function $\hat{G}(\textbf{k},\tau)$ by a quadratic polynomial defined in the interval $[t_{i-1},t_{i+1}]$ 
%%%%%%%%%%%%%%%%%%%%%%%%%%%%%%%%%% 
	\begin{eqnarray}
	\label{G-approx-simpson}
	 \hat{G}(\textbf{k},\tau) &\approx& \hat{G}(\textbf{k},t_{i-1})\frac{(\tau-t_{i})(\tau-t_{i+1})}{2(\Delta t)^2}
	 - \hat{G}(\textbf{k},t_{i})\frac{(\tau-t_{i-1})(\tau-t_{i+1})}{(\Delta t)^2}
	 \nonumber\\
	 &+& \hat{G}(\textbf{k},t_{i+1})\frac{(\tau-t_{i-1})(\tau-t_{i})}{2(\Delta t)^2}.
	\end{eqnarray} 
	%%%%%%%%%%%%%%%%%%%%%%%%%%%%%%%%%%
Substituting Eq.~(\ref{G-approx-simpson}) back into (\ref{I-simpson-int}) and integrating by parts gives a recursive formula for the Duhamel integral \eqref{time-integral}:
	%%%%%%%%%%%%%%%%%%%%%%%%%%%%%%%%%%
	\begin{equation}
	\label{Simpson-recurrence}
	\hat{I}(\textbf{k},t_{i+1}) = e^{2\Delta t \hat{\mathcal{L}}(\textbf{k})}[\hat{I}(\textbf{k},t_{i-1})+q_1\hat{G}(\textbf{k},t_{i-1})+ q_2\hat{G}(\textbf{k},t_{i}) + q_3\hat{G}(\textbf{k},t_{i+1})].	
	\end{equation}
	%%%%%%%%%%%%%%%%%%%%%%%%%%%%%%%%%%
 The quadrature coefficients $q_j\equiv q_j({\bf k},\Delta t), j=1,2,3,$ depend on the Fourier wavenumber and the time 
step $\Delta t$ but not on the iteration index $i$. Thus, they are computed only once.
The exact expressions for the $q_j$'s, $j=1,2,3$, are given by ($z\equiv \Delta t \hat{\mathcal{L}}(\textbf{k})$)
\begin{subequations}\label{quad_coeff}
\begin{eqnarray}
		%\label{D-quad-coeff}
	  \qquad\qquad   q_1 
	   &=&
	 \Delta t(-ze^{-2z}-2e^{-2z}+2z^2-3z+2 )/(2z^3),\\% \nonumber \\
%%%%%%%%%%%%%%%%%%%%%%%%%%%%%%%%%%%%%%%%%%%%%%%%%%%%%%%%%%%%
       %\label{E-quad-coeff}
	    q_2
	    &=&
           \Delta t(2ze^{-2z}+2e^{-2z}+2z-2)/z^3, \\% \nonumber \\ 
%%%%%%%%%%%%%%%%%%%%%%%%%%%%%%%%%%%%%%%%%%%%%%%%%%%%%%%%
	 % \label{F-quad-coeff}
           q_3
           &=&
        \Delta t(-2z^2e^{-2z}-3ze^{-2z}-2e^{-2z}-z+2)/(2z^3). %\nonumber \\ 
%%%%%%%%%%%%%%%%%%%%%%%%%%%%%%%%%%%%%%%%%%%%%%%%%	  
\end{eqnarray}
	   \end{subequations}
	   %%%%%%%%%%%%%%%%%%%%%%%%%%%%%%%%%%%%%%%%%%%%%%%%%
For linear operators satisfying $\hat{\mathcal{L}}(0)=0$ we find (in the limiting case ${\bf k}\rightarrow 0$),
	$q_1(0,\Delta t)  = q_2(0,\Delta t)/4 = q_3(0,\Delta t) \equiv \Delta t/3.$
	   %%%%%%%%%%%%%%%%%%%%%%%%%%%%%%%%%%%%%%%%%%%%%%%%%%
%	$q_4   = q_7 \equiv 3\Delta t/8,$ and 
%	$q_5 = q_6 \equiv 9\Delta t/8.$
	   %%%%%%%%%%%%%%%%%%%%%%%%%%%%%%%%%%%%%%%%%%%%%%%%%%
Equation (\ref{Simpson-recurrence}) needs to be initialized with $\hat{I}(\textbf{k},t=0)=0$ and the quantity $\hat{I}(\textbf{k},\Delta t)=\int_{0}^{\Delta t}e^{(\Delta t-\tau) \hat{\mathcal{L}}(\textbf{k})}\hat{G}(\textbf{k},\tau)d{\tau}$ which we next explain how to find. Note that in the interval $[0,\Delta t]$, the values (in time) of the function $\hat{G}(\textbf{k},\tau)$ are available only at two grid points: $0$ and $\Delta t$. To maintain the same order of accuracy as was done at the other temporal grid points, we apply a combination of two different quadrature rules to approximate $\hat{I}({\bf k},\Delta t)$. First, consider the following identity:
%%%%%%%%%%%%%%%%%%%%%%%%%%%%%%%%%%%%%%%%%%%%%%%%%%%%%%%%%%%%%%%%%%%%%
\begin{equation}
    \label{Calculus-identity-IhatDeltat}
    \int_{0}^{3\Delta t}e^{(\Delta t-\tau) \hat{\mathcal{L}}(\textbf{k})}\hat{G}(\textbf{k},\tau)d{\tau}
    =
     \underbrace{\int_{0}^{\Delta t}e^{(\Delta t-\tau) \hat{\mathcal{L}}(\textbf{k})}\hat{G}(\textbf{k},\tau)d{\tau}}_{\hat{I}(\textbf{k},\Delta t)}
    +
    \int_{\Delta t}^{3\Delta t}e^{(\Delta t-\tau) \hat{\mathcal{L}}(\textbf{k})} \hat{G}(\textbf{k},\tau)d{\tau}.
\end{equation}
%%%%%%%%%%%%%%%%%%%%%%%%%%%%%%%%%%%%%%%%%%%%%%%%%%%%%%%%%%%%%%%%%%%%
The second integral on the right-hand side of Eq.~(\ref{Calculus-identity-IhatDeltat}) is computed using a quadratic interpolation (in time) 
for $\hat{G}(\textbf{k},\tau).$ Indeed, after some algebra, we find
	%%%%%%%%%%%%%%%%%%%%%%%%%%%%%%%%%%
	\begin{equation}
	\label{Filon-simpson-init}
	   \int_{\Delta t}^{3\Delta t}e^{(\Delta t-\tau) \hat{\mathcal{L}}(\textbf{k})}\hat{G}(\textbf{k},\tau)d{\tau} 
	   \approx
	    q_1\hat{G}(\textbf{k},\Delta t)
	   + q_2\hat{G}(\textbf{k},2\Delta t) + q_3\hat{G}(\textbf{k},3\Delta t).
	\end{equation}
	%%%%%%%%%%%%%%%%%%%%%%%%%%%%%%%%%%
To obtain a similar order of accuracy for the integral on the left hand side of Eq.~(\ref{Calculus-identity-IhatDeltat}), we first represent 
$\hat{G}(\textbf{k},\tau)$ as a cubic polynomial (in time) 
	%%%%%%%%%%%%%%%%%%%%%%%%%%%%%%%%%%%%%%%%%
	\begin{eqnarray}
	\label{Filon-simpson-cubic}
	    \hat{G}(\textbf{k},\tau) &\approx&  
	    -\hat{G}(\textbf{k},0)\frac{(\tau-\Delta t)(\tau-2\Delta t)(\tau-3\Delta t)}{6(\Delta t)^3}
	    + \hat{G}(\textbf{k},\Delta t)\frac{\tau (\tau-2\Delta t)(\tau-3\Delta t)}{2(\Delta t)^3} 
	    \nonumber \\
	    &-&
	   \hat{G}(\textbf{k},2\Delta t)\frac{\tau (\tau-\Delta t)(\tau-3\Delta t)}{2(\Delta t)^3}
	    +
	   \hat{G}(\textbf{k},3\Delta t)\frac{\tau (\tau-\Delta t)(\tau-2\Delta t)}{6(\Delta t)^3} \;.
	\end{eqnarray}
	%%%%%%%%%%%%%%%%%%%%%%%%%%%%%%%%%%%%%%%%%
 Substituting expressions (\ref{Filon-simpson-cubic}) and (\ref{Filon-simpson-init}) into Eq.~(\ref{Calculus-identity-IhatDeltat}) gives (after integration by parts)
	%%%%%%%%%%%%%%%%%%%%%%%%%%%%%%%%%%%%%%%%%
	\begin{eqnarray}
	\label{Exact-formula-initialize}
	    \hat{I}(\textbf{k},\Delta t) 
	    &=&  
	    q_4 e^{\Delta t \hat{\mathcal{L}}(\textbf{k})} \hat{G}(\textbf{k},0)
	    +
	     \left( q_5 e^{\Delta t \hat{\mathcal{L}}(\textbf{k})} 
	     - q_1\right)\hat{G}(\textbf{k},\Delta t)
	    \\ \nonumber
	    &+&\left( q_6 e^{\Delta t \hat{\mathcal{L}}(\textbf{k})} - q_2\right)\hat{G}(\textbf{k},2\Delta t)
	    +\left(q_7 e^{\Delta t \hat{\mathcal{L}}(\textbf{k})} - q_3\right)\hat{G}(\textbf{k},3\Delta t).	
	    \end{eqnarray}
%%%%%%%%%%%%%%%%%%%%%%%%%%%%%%%%%%%%%%%%%%%%%%%%%%%%%%
Here, $q_j\equiv q_j({\bf k},\Delta t), j=4,5,6,7$, denote the quadrature coefficients whose expressions are given by
\begin{subequations}
\begin{eqnarray}
     % \label{H-quad-coeff}
	   q_4 &=&
	   \Delta t(2z^2e^{-3z}+6ze^{-3z}+6e^{-3z}+6z^3+12z-11z^2-6)/(6z^4), \\% \nonumber \\
%%%%%%%%%%%%%%%%%%%%%%%%%%%%%%%%%%%%%%%%%%%%%%%
	 %  \label{M-quad-coeff}
           q_5
            &=& 
           \Delta t (-3z^2e^{-3z}-8ze^{-3z}-6e^{-3z}+6z^2-10z+6)/(2z^4),  \\%\nonumber \\ 
	  %%%%%%%%%%%%%%%%%%%%%%%%%%%%%%%%%%%%%%
	   %\label{N-quad-coeff}
	   q_6 
	   &=& 
	   \Delta t(6z^2e^{-3z}+10ze^{-3z}+6e^{-3z}-3z^2+8z-6)/(2z^4),  \\ %\nonumber \\ 
	%%%%%%%%%%%%%%%%%%%%%%%%%%%%%%%%%%%%%%%%%%%
	%  \label{O-quad-coeff}
	  q_7 
	  &=&
	  \Delta t(-6z^3e^{-3z}-11z^2e^{-3z}-12ze^{-3z}-6e^{-3z}+2z^2-6z+6)/(6z^4).
\end{eqnarray}
\label{quad_coeff-2}
\end{subequations}
For linear operators satisfying $\hat{\mathcal{L}}(0)=0$, and for  wavenumber ${\bf k}\rightarrow 0$, $q_4(0,\Delta t)   = q_7(0,\Delta t) \equiv 3\Delta t/8,$ and 
	$q_5(0,\Delta t) = q_6(0,\Delta t) \equiv 9\Delta t/8.$
	   %%%%%%%%%%%%%%%%%%%%%%%%%%%%%%%%%%%%%%%%%%%%%%%%%%
To summarize, the Duhamel integral $I({\bf x},t)$ is determined from iterating Eq.~(\ref{Simpson-recurrence}) subject to the initial conditions: $I(\textbf{x},0)=0$ and $I(\textbf{x},\Delta t)$ given in Fourier space by Eq.~(\ref{Exact-formula-initialize}).  Note that in some cases, the Filon coefficients {$q_j$} may exhibit a removable singularity in the variable $z\equiv \Delta t \hat{\mathcal{L}}({\bf k})$ at zero wave number that could trigger numerical instability. To remedy this, we represent each quadrature term as a Cauchy integral that allows a stable and uniform approximation valid for all wavenumbers. This idea has been first implemented in the context of exponential time differencing fourth order Runge-Kutta (ETDRK4) \cite{kassam2005fourth}.
 %%%%%%%%%%%%%%%%%%%%%%%%%%%%%%%%%%%%%%%%%%%%%%%%
For the sake of completeness, we show how to implement this approach on the coefficient $q_1.$ The computation of the other quadrature coefficients follow similar derivation. 
%%%%%%%%%%%%%%%%%%%%%%%%%%%%%
Since the function ${q_1}(\zeta;\Delta t)$ is analytic in the $\zeta$ complex plane, by the Cauchy integral formula we have
%%%%%%%%%%%%%%%%%%%%%%%%%%%%%
\begin{equation}
  \label{D-quad-cauchy}
    q_1(z;\Delta t)= \frac{1}{2\pi i} \int_{\mathcal{C}} \frac{q_1(\zeta;\Delta t)}{\zeta-z}d{\zeta},
\end{equation}
%%%%%%%%%%%%%%%%%%%%%%%%%%%%%
where $\mathcal{C}$ is a circle of constant radius centered at $z$. The above integral can be evaluated to spectral accuracy with the use of the trapezoidal quadrature  \cite{kassam2005fourth,davis1959numerical}.
%%%%%%%%%%%%%%%%%%%%%%%%%%%%%%
%Having formulated the TDSR algorithm subject to periodic and decaying boundary conditions, we next direct our attention to the numerical implementation. 

%%%%%%%%%%%%%%%%%%%%%%%%%%%%%%%%%%%%%%%%%%%%%%%%%%%%%%%
%\st{Numerical tests indicate that the algorithm converges to the correct solution when using the root \\
 %$R_2(t) = [-\mu_2(t)+\sqrt{\mu_2^2(t)-4\mu_1(t) \mu_3(t)}]/[2 \mu_1(t)]$} while the other root results to \scn{divergence of the TDSR algorithm}.\scn{The ``correct" root was found to satisfy $R_2(0)v_2(x,0)=f_2(x)$ at every Duhamel iteration; the incorrect root consistently violated this in all our numerical tests.} It \scn{was seen} that the mass and momentum are indeed conserved by the TDSR. \st{This is in contrast with other well established methods where conservation of a single quantity only is possible.}\scn{NOTE: THIS IS NOT A TRUE STATEMENT; THERE ARE SOME SOME CITED METHODS WHICH CONSERVE THE MASS AND MOMENTUM FOR THE KdV} 
%%%%%%%%%%%%%%%%%%%%%%%%%%%%%%%%%%%%%%%%%%%%%
%%%%%%%%%%%%%%%%%%%%%%%%%%%%%%%%%%%%%%%%%%%%%%%%%%%%%%%%%%%
\subsection{Time integration: non-periodic boundary conditions}
\label{time-integration-non-periodic}
 %%%%%%%%%%%%%%%%%%%%%%%%%%%%%%%%%%%%%%%%%%%%%%%%%%%%%%%%%%%
So far, we have discussed the development and application of the TDSR method to evolution equations subject to periodic or localized boundary conditions. Here, we intend to extend the TDSR scheme to allow for non-periodic and non-decaying boundary conditions where the use of Fourier analysis is not applicable. In such circumstances the matrix approximating the linear operator could be banded (as is the case with finite differences) or dense, for example, in case of Chebyshev spectral method. 
%While the previous example of KdV dispersive shock waves falls into this category, it was nonetheless solved using TDSR subject to decaying boundary conditions after taking its spatial derivative. 
The derivation of the Duhamel formula follows similar steps as outlined in Sec.\ref{Sec-TDSR} with the exception of the use of trapezoidal scheme instead of Simpson. 
%The integral we are interested in is given in 
%Eq.~(\ref{time-integral-non-periodic-exact}) which, for the ease of presentation, we rewrite it again:
%      \begin{equation}
%    \label{time-integral-non-periodic-exact}
%        I=\int_{0}^{t}e^{(t-\tau) \mathcal{L}} G({x},\tau)d{\tau}.
%    \end{equation}
%%%%%%%%%%%%%%%%%%%%%%%%%%%%%%%%%%%%%%%%%%%%%%%%%%%%%%%%
%To approximate the above integral in Eq.~(\ref{time-integral}), we start by discretizing the space variable $x\rightarrow x_j$ and $G(x,\tau)\rightarrow {\bf G}(\tau)\equiv G(x_j,\tau), \; j=1, 2, \cdots, N_S$. 
Using a Chebyshev basis function or other discretization methods we represent the differential operator $\mathcal{L}$ in Eq. \eqref{time-integral} by a finite dimensional matrix ${\bf L}$. The boundary conditions are incorporated within the matrix ${\bf L}$. 
%With this at hand, we seek to approximate the Duhamel time integral
%    \begin{equation}
%    \label{time-integral-non-periodic}
%        {\bf I}(t)=\int_{0}^{t}e^{(t-\tau) {\bf L}} {\bf G}(\tau)d\tau.
%   \end{equation}
%%%%%%%%%%%%%%%%%%%%%%%%%%%%%%%%%%%%%%%%%%%%%%
By creating a mesh in time domain, Eq.\eqref{time-integral} can be put in the recursive form 
%%%%%%%%%%%%%%%%%%%
\begin{equation}
\label{Recurrence-exact-nondiagonal}
    {\bf I}(t_{i+1})= e^{\Delta t {\bf L}} {\bf I}(t_i) + e^{\Delta t {\bf L}} \int_{t_i}^{t_{i+1}} e^{(t_{i}-\tau) {\bf L}} {\bf G}(\tau) d{\tau},
\end{equation}
%%%%%%%%%%%%%%%%%%%%%%%%%%%%%%%%%%%%%%%%%%%%%%%
where now ${\bf I}(t_i)$ is the matrix representing the Duhamel integral at space meshgrid $x$ and time level $t_i$. Additionally, ${\bf G}(t_i)$ is the matrix representing the nonlinear terms at $t_i=i\Delta t$ with $i=0,1,\cdots N_T$. Using a linear interpolant to approximation ${\bf G}(\tau)$ in the interval $[t_i, t_{i+1}]$ we find
%%%%%%%%%%%%%%%%%%%%%%%%%%%%%%%%%%%%%%%%%%%%%%%%
\begin{equation}
\label{linear-interpolant-non-diagonal}
{\bf G}(\tau)\approx {\bf G}(t_i)+\frac{{\bf G}(t_{i+1})-{\bf G}(t_i)}{\Delta t}(\tau-t_i).
\end{equation}
%%%%%%%%%%%%%%%%%%%%%%%%%%%%%%%%%%%%%%%%%%%%%%%%%
Substituting Eq.(\ref{linear-interpolant-non-diagonal}) into (\ref{Recurrence-exact-nondiagonal}) we obtain after some algebra
%%%%%%%%%%%%%%%%%
\begin{equation}
\label{}
 {\bf I}(t_{i+1})= e^{\Delta t {\bf L}}[{\bf I}(t_i) + A {\bf G}(t_i) + B {\bf G}(t_{i+1})],
\end{equation}
%%%%%%%%%%%%%%%%%%
where the matrix valued quadrature coefficients $A\equiv A({\bf L},\Delta t)$ and $B\equiv B({\bf L},\Delta t)$ are defined by
%%%%%%%%%%%%%%%%%%%%%%%%%%%%%%%%%%%%%%%%%%%%%%%%%%%%
\begin{align}
    A &\equiv \Delta t {\tilde A},\quad  \tilde{A}={\bf \tilde{L}}^{-2}\left(e^{- {\bf \tilde{L}}} + {\bf \tilde{L}} - \mathcal{I}\right), \label{A-quad-matrix-coeff} \\
    B &\equiv \Delta t {\tilde B}, \quad \tilde{B}={\bf \tilde{L}}^{-2}\left( \mathcal{I} - {\bf \tilde{L}} e^{-{\bf \tilde{L}}}-e^{-{\bf \tilde{L}}}\right), \label{B-quad-matrix-coeff}
\end{align}
with $\mathcal{I}$ being the identity matrix and $\tilde{{\bf L}} \equiv \Delta t {\bf L}$. 
As was done in Sec.~\ref{time-int-periodic} for the periodic case \cite{kassam2005fourth}, we again adopt the Cauchy integral formula to represent each quadrature coefficient as a contour integral in the complex plane. Thus we write :
%%%%%%%%%%%%%%%%%%%%%%%%%%%%%%%%%%%%%%%%%%%%%%%%%%%%%%%%%%%%%%
%%%%%%%%%%%%%%%%%%%%%%%%%%%%%%%%%%%%%%%%%%%%%%%%%%%%%%
\begin{equation}
\label{cauchy-filon-B-matrix}
\tilde{A}( \tilde{{\bf L}} ) = \frac{1}{2\pi i}\int_{\Gamma} \tilde{A}(\zeta) (\zeta\mathcal{I} - \tilde{{\bf L}})^{-1} d{\zeta},\;\;\;\tilde{B}( \tilde{{\bf L}} ) = \frac{1}{2\pi i}\int_{\Gamma} \tilde{B}(\zeta) (\zeta\mathcal{I} - \tilde{{\bf L}})^{-1} d{\zeta}
\end{equation}
%%%%%%%%%%%%%%%%%%%%%%%%%%%%%%%%%%%%%%%%%%%%%%%%%%%%%%%
with $\Gamma$ being a circular contour that encloses all the eigenvalues of $\tilde{{\bf L}}$. The integral in Eq.(\ref{cauchy-filon-B-matrix}) is computed to a spectral accuracy with the use of the trapezoidal rule. 

\bigskip
\bigskip
%%%%%%%%%%%%%%%%%%%%%%%%%%%%%%%%%%%%%%%%%%%%%%%%%%%%
\section{Numerical Implementation of TDSR: One conservation or dissipation law with various boundary conditions}
%%%%%%%%%%%%%%%%%%%%%%%%%%%%%%%%%%%%%%%%%%%%%%%%%%%%

%%%%%%%%%%%%%%%%%%%%%%%%%%%%%%
\subsection{The KdV equation.}
\begin{figure}
\begin{subfigure}%{0.33\textwidth}
  \centering
  \includegraphics[width=50mm]{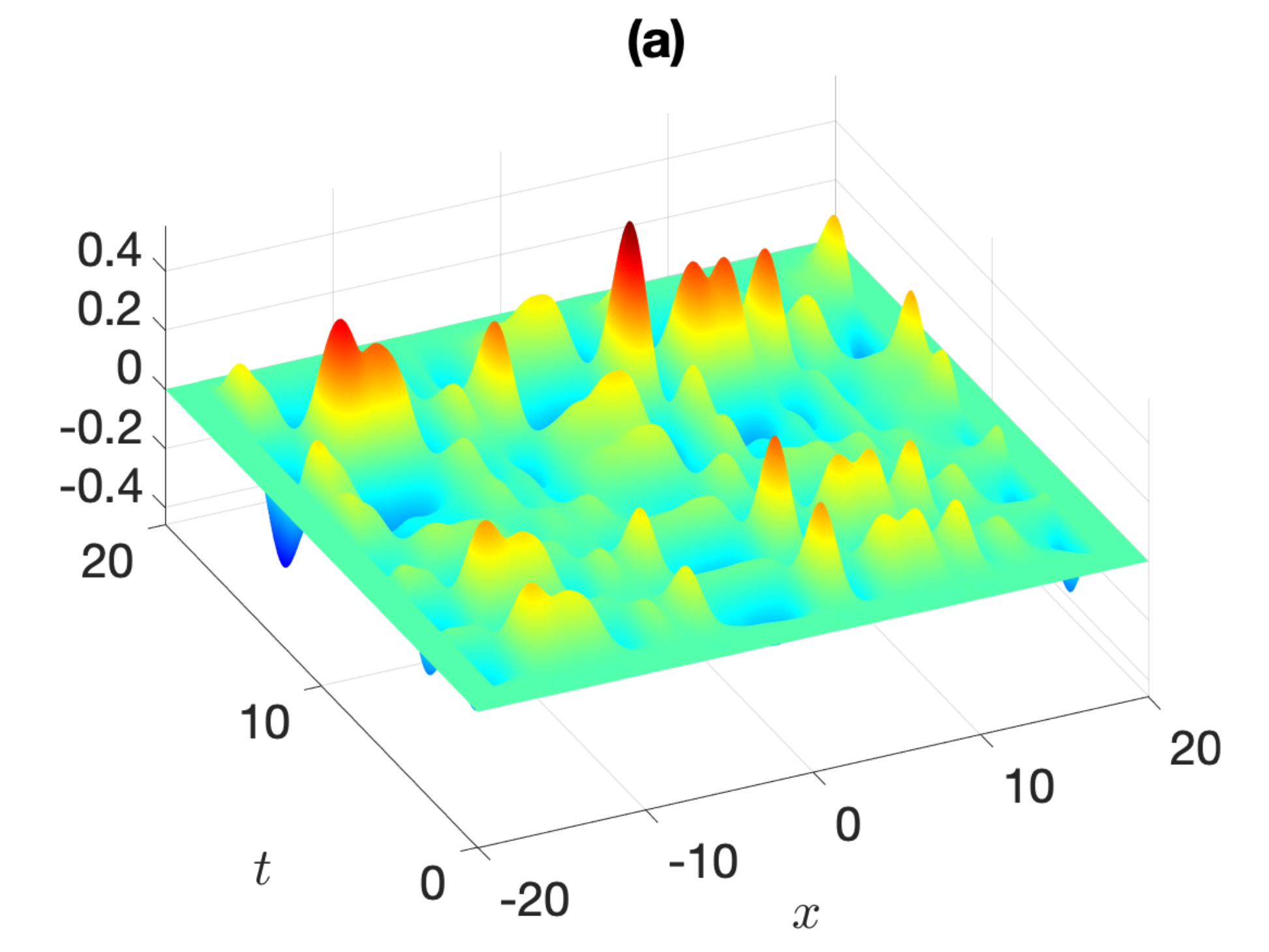}
%  \caption{Initial guess: random space-time function modulated by a bump function of spatial extent $l=36$}
  %\label{fig:sfig1}
\end{subfigure}%
\begin{subfigure}%{0.34\textwidth}
  \centering
  \includegraphics[width=55mm]{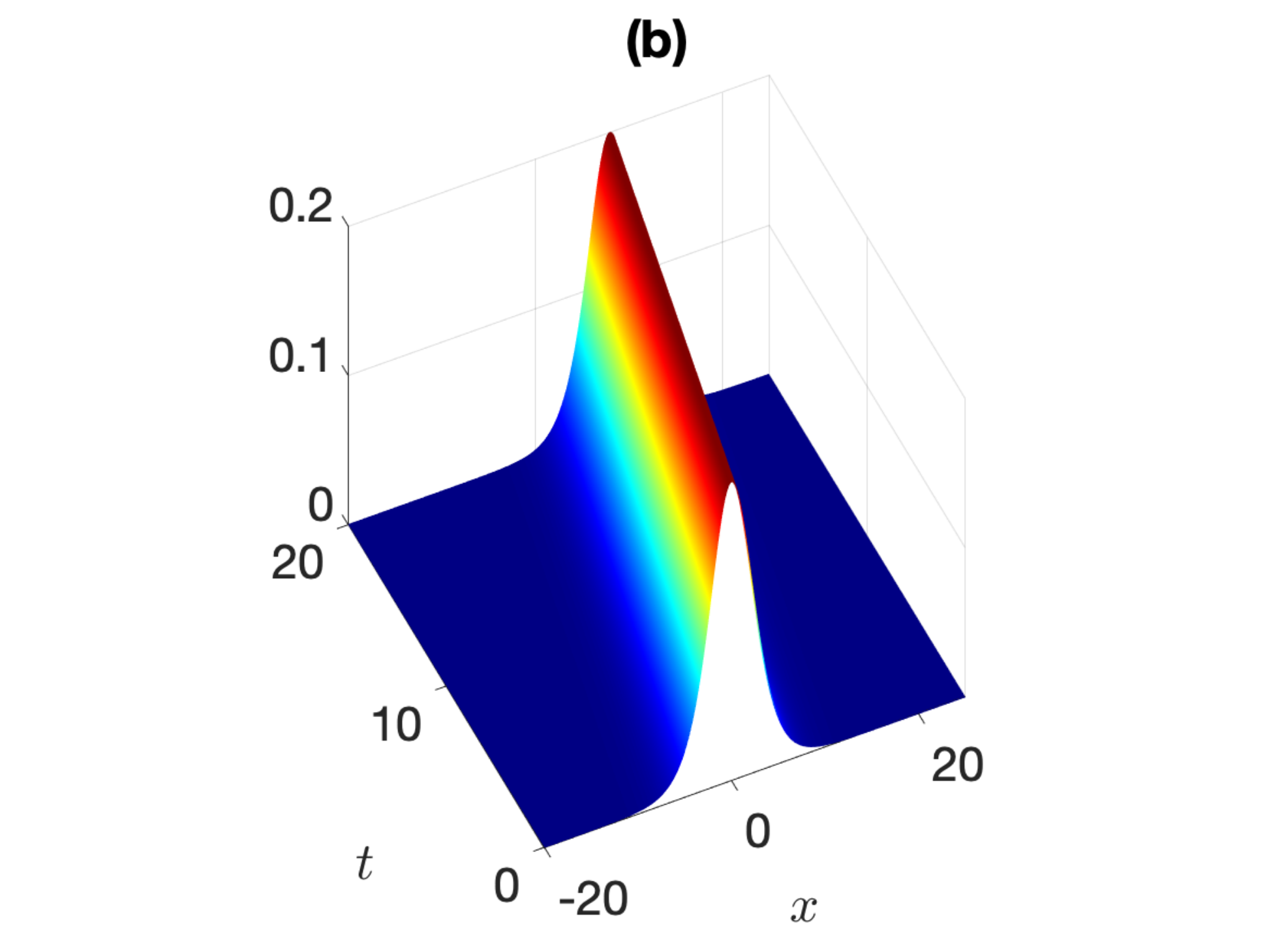}
 % \caption{Results upon convergence after 30 Duhamel iterations}
  %\label{fig:sfig2}
\end{subfigure}
\begin{subfigure}%{0.32\textwidth}
  \centering
  \includegraphics[width=35mm]{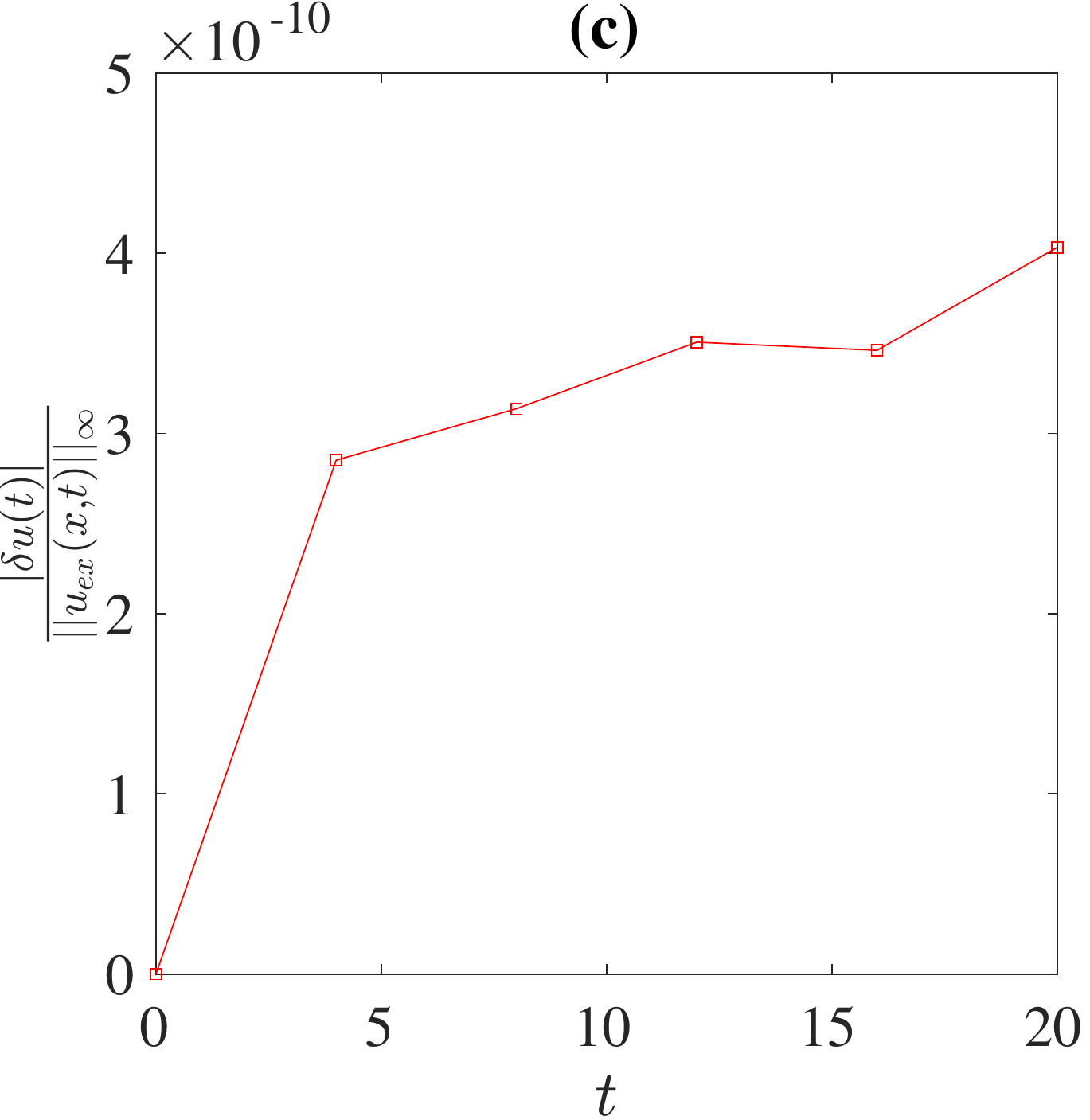}
 % \caption{Results upon convergence after 30 Duhamel iterations}
  %\label{fig:sfig2}
\end{subfigure}
%\vspace{-0.1cm}
\caption{\small
(a) A random space-time initial guess constructed from a linear superposition of randomly centered Gaussians with random amplitudes -- see 
Eq.(\ref{Space-time-arbitrary-ig}). (b) Numerical solution for the KdV equation obtained from the TDSR algorithm after 30 Duhamel iterations. Parameters are: $T=20, \Delta t = 0.025, L=100, N_S=2048$ (Fourier modes). Here, the wave speed is $4\beta^2 = 2/5$. This figure was generated by imposing conservation of momentum for which the renormalization parameter $R(t)$ is computed from Eq.~(\ref{Renormalization-factor-momentum}). The soliton initial condition is $u_0(x) = 2\beta^2{\rm sech}^2(\beta x)$ with $\alpha=6$ and $\epsilon =1.$ (c) Time evolution of the relative error between the TDSR and the exact solution. Mollifier parameters are $a=0.95\times\frac{L}{2}$, $b=1.$}
\label{1-soliton-simp-result}
\end{figure}

\begin{figure}
\begin{subfigure}%{0.49\linewidth}
 \centering
\includegraphics[width=75mm]{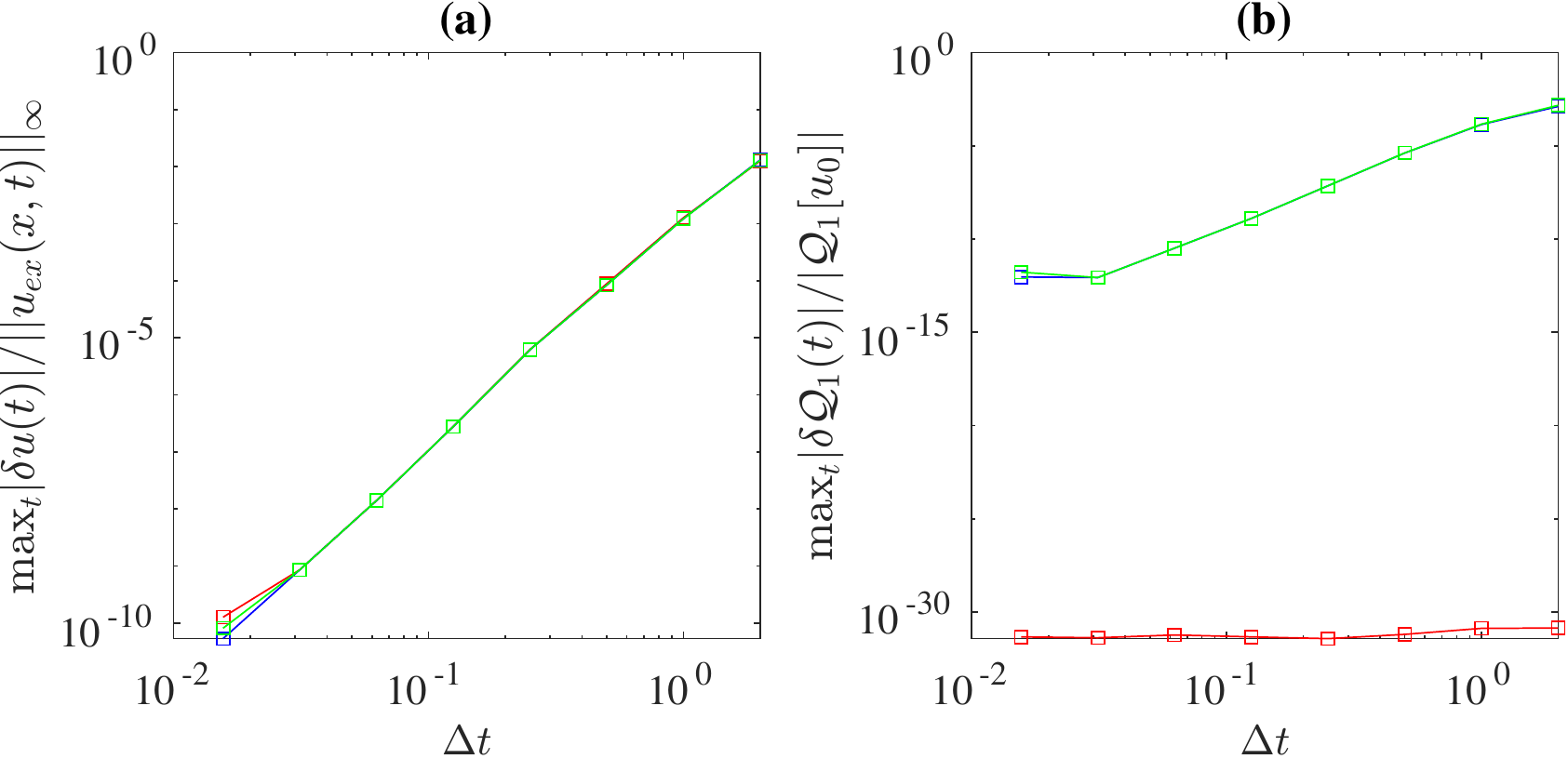}
%\caption{(a) Infinity norm of error in numerical solution (b) Error in }
\end{subfigure}
\begin{subfigure}%{0.49\linewidth}
 \centering
\includegraphics[width=75mm]{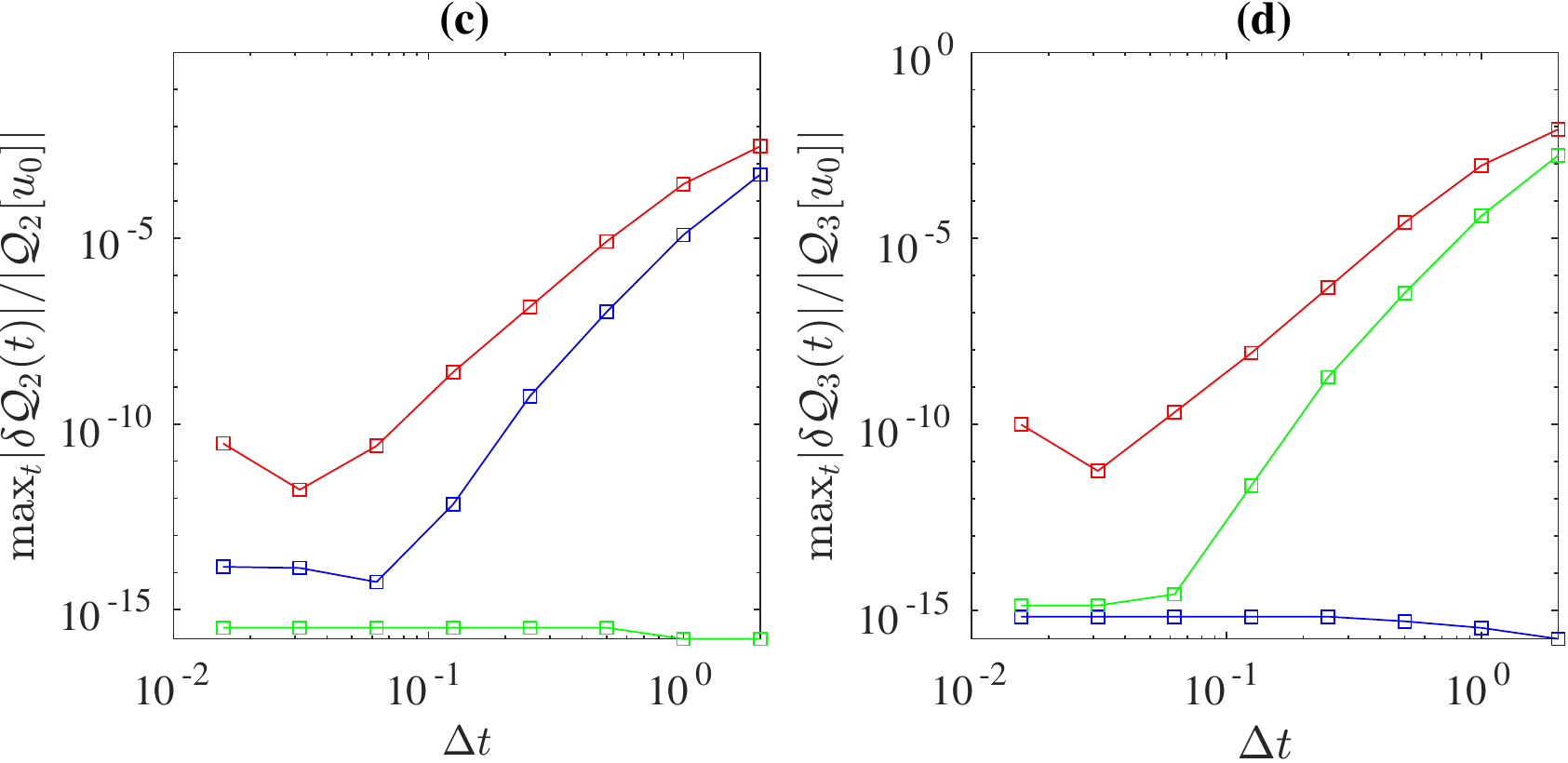}
%\caption{Image B}
\end{subfigure}%
%\vspace{-0.4cm}
\caption{\small The relative error in (a) TDSR solution, (b) mass, (c) momentum, and (d) Hamiltonian. Parameters are $T=10$, $L=100$, $N_S=2048$. The renormalization factor $R(t)$ is computed by enforcing either conservation of mass (red), momentum (green), or Hamiltonian (blue).}
\label{TDSR-Simpson-family}
\end{figure}

%%%%%%%%%%%%%%%%%%%%%%%%%%%%%%
In this section, we use the KdV equation as a testbed PDE model to examine various numerical aspects related to the TDSR method such as convergence, accuracy and dependence on initial guesses. 
The KdV equation is given by:
\begin{eqnarray}
\label{kdv-eqn}
u_t + \alpha uu_x + \epsilon^2 u_{xxx}=0,
\end{eqnarray}
where $\alpha, \epsilon$ are real and positive numbers. When considered on the whole real line with rapidly decaying boundary conditions, Eq.(\ref{kdv-eqn}) admits a one parameter family of soliton solution given by (e.g. $\alpha=6, \epsilon =1$)
\begin{equation}
\label{one-para-soliton}
u_{ex}(x,t)= 2 \beta^2 {\rm sech}^2(\beta (x-4\beta^2t)), \;\;\; \beta >0.
\end{equation}
It is noteworthy that the KdV equation is an integrable dynamical system admitting infinitely many conservation laws. Among them are the physically relevant mass, momentum and Hamiltonian given for $\alpha=6, \epsilon =1$ by 
Eq.~(\ref{n-cons-laws}) with $Q_1 = u,\; Q_2 = u^2$ and $Q_3= -u^3 + \frac 12u_x^2$, respectively.
All numerical simulations reported in this section were performed on a spatial domain of size $L=100$ or $L=800$ ({depending on the case at hand}) with corresponding number of spatial grid points (Fourier modes) $N_S=2048$, $N_S=16384$ respectively and time interval  $[0, T]$ with $T = 5, 10, 20, 30$ or $60$. In this section, the renormalization factors were computed by enforcing a single conservation law. Numerical convergence and accuracy were quantified by monitoring the error between two successive iterations 
$$\underset{x,t}{\rm max}|u^{(n+1)}(x,t)-u^{(n)}(x,t)|,$$
and the quantities
%%%%%%%%%%%%%%%%%%%%%%%%%%%%%%%
\begin{align}
\delta u(t) \equiv \underset{x}{{\rm max}} |u(x,t)-u_{ex}(x,t)|,\label{solution-accuracy} \\
\delta \mathcal{Q}_j(t) \equiv \mathcal{Q}_j[u(x,t)] - \mathcal{Q}_j[u_0(x)] \;,\;\;\; j=1,2,3, \label{cons-law-convergence}
\end{align}
where $u_0(x)=2 \beta^2 {\rm sech}^2(\beta x)$ is the initial condition associated with Eq.~(\ref{kdv-eqn}), for the parameters
$\alpha=6, \epsilon =1$, {where $u_{ex}(x,t)$ is the one-parameter soliton solution for the KdV given in Eq.~(\ref{one-para-soliton})}. For all simulations reported here, convergence tolerance was set near $1\times10^{-16}$.
All functionals $\mathcal{Q}_j, j=1,2,3,$ are defined in Eq.~(\ref{n-cons-laws}).
 %%%%%%%%%%%%%%%%%%%%%%%%%%%%%%%%%%%%%%%%%%%%%%%%%%%%%%%%%%%
We perform numerical experiments on the KdV equation while conserving one of the following quantities: mass, momentum, or Hamiltonian. The renormalization factor $R(t)$ in each case, is given by: 
%%%%%%%%%%%%%%%%%%%%%%%%%%%%%%%%%%%%%%%%
\begin{eqnarray}
    \label{Renormalization-factor-mass}
   && \text{mass:}
    \;\;\;\;\;\;\;\;\;\;\;\;\;\;
     R(t) = \frac{\mathcal{Q}_1[u_0(x)]}{\mathcal{Q}_1[v(x,t)]},
     \\ 
    \label{Renormalization-factor-momentum}
    &&\text{momentum:}  \;\;\;\;
     R(t) = \left(\frac{\mathcal{Q}_2[u_0(x)]}{\mathcal{Q}_2[v(x,t)]}\right)^{1/2},
     \\
    \label{Renormalization-factor-Hamiltonian}
    &&\text{Hamiltonian:}
    \;\;\;
     A(t)R^3(t)-B(t)R^2(t)-C_3=0,
    \end{eqnarray}
where $A(t) = -\int_{\mathbb{R}} v^3(x,t)\,dx$, $B(t) =  -\frac{1}{2} \int_{\mathbb{R}} v_x^2\,dx$ and $C_3$ is the initial value of the Hamiltonian.
All spatial integrals are computed to spectral accuracy with the use of fast Fourier transform.
We initialize the TDSR algorithm with a space-time random function $v^{(1)}({ x},t)$ constructed by superimposing several Gaussians each being centered at a random location and having a random time-dependent amplitude. The centers and amplitudes are sampled from a uniform distribution on the interval $[-L/2, L/2]$ and $[-1, 1]$ respectively. To ensure the initial guess $v^{(1)}({x}, t)$ satisfies the underlying boundary conditions, we mollify it with $\chi (x)$. Thus, we have 
%%%%%%%%%%%%%%%%%%%%%%%%%%%%%%%%%%
\begin{equation}
\label{Space-time-arbitrary-ig}
    {v}^{(1)}(x,t)
    =
    \frac{\sum_{n=1}^{N_G}a_{n}(t)\exp\big[-\big(\frac{x-c_{n}}{d}\big)^{2}\big]}{{\rm max}_{x,t} |\sum_{n=1}^{N_G}a_{n}(t)\exp\big[-\big(\frac{x-c_{n}}{d}\big)^{2}\big]|}
     \chi (x),
\end{equation}
where $N_G$ denotes the number of Gaussians with $c_n$ and $a_n(t)$ representing their centres and time-varying amplitudes and $d$ is the width. Here, $\chi(x)$ is the mollifier with unit peak amplitude defined by:
\begin{equation}
\label{Mollifier}
\chi(x) = \begin{cases}
 {\rm exp}\bigg({\frac{b}{a^2} -\frac{b}{a^2-x^2}}\bigg), {\rm if}\hspace{1mm} x\in (-a,a)\\
 0, {\rm if} \hspace{1mm}|x|>a,
\end{cases}
\end{equation}
with arbitrary mollifier parameters $a$ and $b$. 
As expected, the numerical result agrees well with the exact solution. In generating Fig.~\ref{1-soliton-simp-result}, conservation of momentum is imposed in which case the renormalization parameter $R(t)$ is computed from Eq.~(\ref{Renormalization-factor-momentum}).
One could instead reach the same conclusion by using a different dynamic renormalization process emanating from either conservation of mass or Hamiltonian.
%%%%%%%%%%%%%%%%%%%%%%%%%%%%%%%%%%%%%%%%%%%%%%%%%%%%%%%%%%%
%%%%%%%%%%%%%%%%%%%%%%%%%%%%%%%%%%%%%%%%%%%%%%%%%%%%%%%%%%%
%%%%%%%%%%%%%%%%%%%%%%%%%%%%%%%%%%%%%%%%%%%%%%%%%%%%%%%
This numerical experiment reveals the simplistic (albeit powerful) nature of our proposed method as measured by its easy formulation, actual implementation, robustness to initial guesses and its ability to impose conservation laws ``on-demand".
To further characterize the numerical performance of the TDSR scheme, we have investigated its temporal convergence properties by measuring the space-time maximum error between the numerically obtained solution to the KdV equation (relative to its exact solution) and all conservation laws, as quantified by Eqs.~(\ref{solution-accuracy}) and (\ref{cons-law-convergence}), as a decreasing function of time step $\Delta t$. It is evident from Fig.~\ref{TDSR-Simpson-family}(a) that the 
maximum (over time) error in the solution decreases at a fourth order rate. This trend seems to persist independently of the choice of specific conservation law. When conservation of mass is imposed, the error in the Hamiltonian and momentum reduces with decreasing $\Delta t$. Similar scenarios occur when imposing conservation of momentum or Hamiltonian -- see Fig.~\ref{TDSR-Simpson-family} (b)-(d).
%%%%%%%%%%%%%%%%%%%%%%%%%%%%%%%%%%%%%%%%%%%%%%%%%%%%%%%%%
%%%%%%%%%%%%%%%%%%%%%%%%%%%%%%%%%%%%%%%%%%%%%%%%%%%%%%%%%%
We remark that to conserve the Hamiltonian structure, we need to solve a cubic equation defined by Eq.~(\ref{Renormalization-factor-Hamiltonian}). 
As such, there are three possible expressions for the renormalization factor, of which only one is feasible. It turns out that the right expression yields $R(0) v(x,0)=u_0(x)$ at any Duhamel iteration, while the other two roots violate this criterion.
Finally, for cases where the renormalization factor satisfies an associated equation that lacks exact solution, one needs to resort to a root finding algorithms such as the Newton's method. It is interesting to note that when it comes to long time simulations, the TDSR performs optimally when imposing conservation of momentum rather than mass or Hamiltonian. Indeed, we have tested the TDSR method on the long-time evolution of the 1-soliton solution for the KdV equation while conserving momentum (the $L^2$ norm of the solution). The numerical experiment was performed with parameters $T=240, \Delta t=0.1875, N_S=4096, L=300$  using the idea of {multi-blocking} with $M_b=8$ time blocks (see remark below). The relative error in the solution, mass and Hamiltonian at end time were in the order of $10^{-6}$, $10^{-7}$ and $10^{-10}$ respectively, while the relative error in momentum remained near machine precision.\\
{\bf Remark}: Below, we describe the idea of multi-blocking used when the time interval is too large for the renormalized Picard iterations to converge (this is not due to a CFL-type restriction prevalent in generic time-stepping schemes). The idea is to divide the full time interval $[0,T]$ into $M_b$ sub-intervals such that $[0, T] = \cup_{i=1}^{M_b} [T_{i-1}, T_i]$ with $T_0=0.$ For a fixed $i$, the quantity $T_{i} - T_{i-1}$ is chosen sufficiently large so that the spectral
renormalization algorithm is efficient and convergent.{ On the first segment $[0,T_1]$, the solution of the TDSR scheme with the initial condition $u_0(x)$ is obtained from the iteration: 
\begin{align*}
    &v^{(n+1)}(x,t)=\frac{1}{R^{(n)}(t)}\Big(e^{t\mathcal{L}}[u_0(x)]+\int_{0}^{t} e^{(t-\tau)\mathcal{L}}\mathcal{N}[R^{(n)}(\tau)v^{(n)}(x,\tau)] d{\tau}\Big),\; \\\nonumber
    %v^{(n+1)}(x,t)&=\frac{1}{R(t)}\Big(e^{(t-T_1)\mathcal{L}}[u(x,T_1)]+\int_{T_1}^{t} e^{(t-T_1-\tau)\mathcal{L}}\mathcal{N}[R(\tau)v(x,\tau)] d{\tau}\Big),\;t\in [T_1,T_2], 
\end{align*}
while on the second segment $[T_1,T_2]$, from: 
\begin{align*}
    &v^{(n+1)}(x,t)=\frac{1}{R^{(n)}(t)}\Big(e^{(t-T_1)\mathcal{L}}[u(x,T_1)]+\int_{T_1}^{t} e^{(t-T_1-\tau)\mathcal{L}}\mathcal{N}[R^{(n)}(\tau)v^{(n)}(x,\tau)] d{\tau}\Big).\; 
\end{align*} 
The renormalization factor $R^{(n)}(t)$ (corresponding to a single conservation law) is obtained from 
\begin{align*}
    &\mathcal{Q}_m\Big(R^{(n)}(t)v^{(n)}(x,t)\Big)=C_m.
    %v^{(n+1)}(x,t)&=\frac{1}{R(t)}\Big(e^{(t-T_1)\mathcal{L}}[u(x,T_1)]+\int_{T_1}^{t} e^{(t-T_1-\tau)\mathcal{L}}\mathcal{N}[R(\tau)v(x,\tau)] d{\tau}\Big),\;t\in [T_1,T_2]
\end{align*}
}

This process is repeated $M_b$ times until final time $T$ is reached. It should be pointed out that the number of segment $M_b$ is chosen such that the Duhamel fixed point iteration, without renormalization, would diverge on any given sub interval $[T_{i-1}, T_i]$.

%\begin{figure}
%    \centering
%    \includegraphics[width=110mm]{Ham_renorm_factors.pdf}
%    \caption{\small The figure summarizes the discussion on the evaluation of the renormalization factor from Eq.\eqref{Renormalization-factor-Hamiltonian}, when the TDSR was initialized with a space-time initial guess $v^{(1)}(x,t)=f(x)$ and computational parameters $T=20$,$\Delta t=0.125$, $N_S=2048$ and $L=100$ (a) The  ``correct" renormalization factor upon algorithm convergence (in red) is indistinguishable from the values computed from Newton's method (overlaid blue squares), the other two renormalization factors were found to be complex valued functions $\tilde R(t)$ and $\tilde R^{*}(t)$. Here $*$ denotes the complex conjugate of a complex number. The successive two figures in the panel are (b) the real and (c) imaginary parts of $\tilde R(t)$ plotted as a function of time.}
%    \label{Ham-renorm-factors}
%\end{figure}
%
%%%%%%%%%%%%%%%%%%%%%%%%%%%%%%%%%%%%%%%%%%%%
\bigskip
%%%%%%%%%%%%%%%%%%%%%%%%%%%%%%%%%%%%%%%%%%%%%%
\subsection{Zabusky-Kruskal experiment }
\begin{figure}%[!hb]
    \begin{subfigure}%[h]{0.49\linewidth}
    \centering
    %\begin{minipage}{0.6\textwidth}
    \includegraphics[width=45mm]{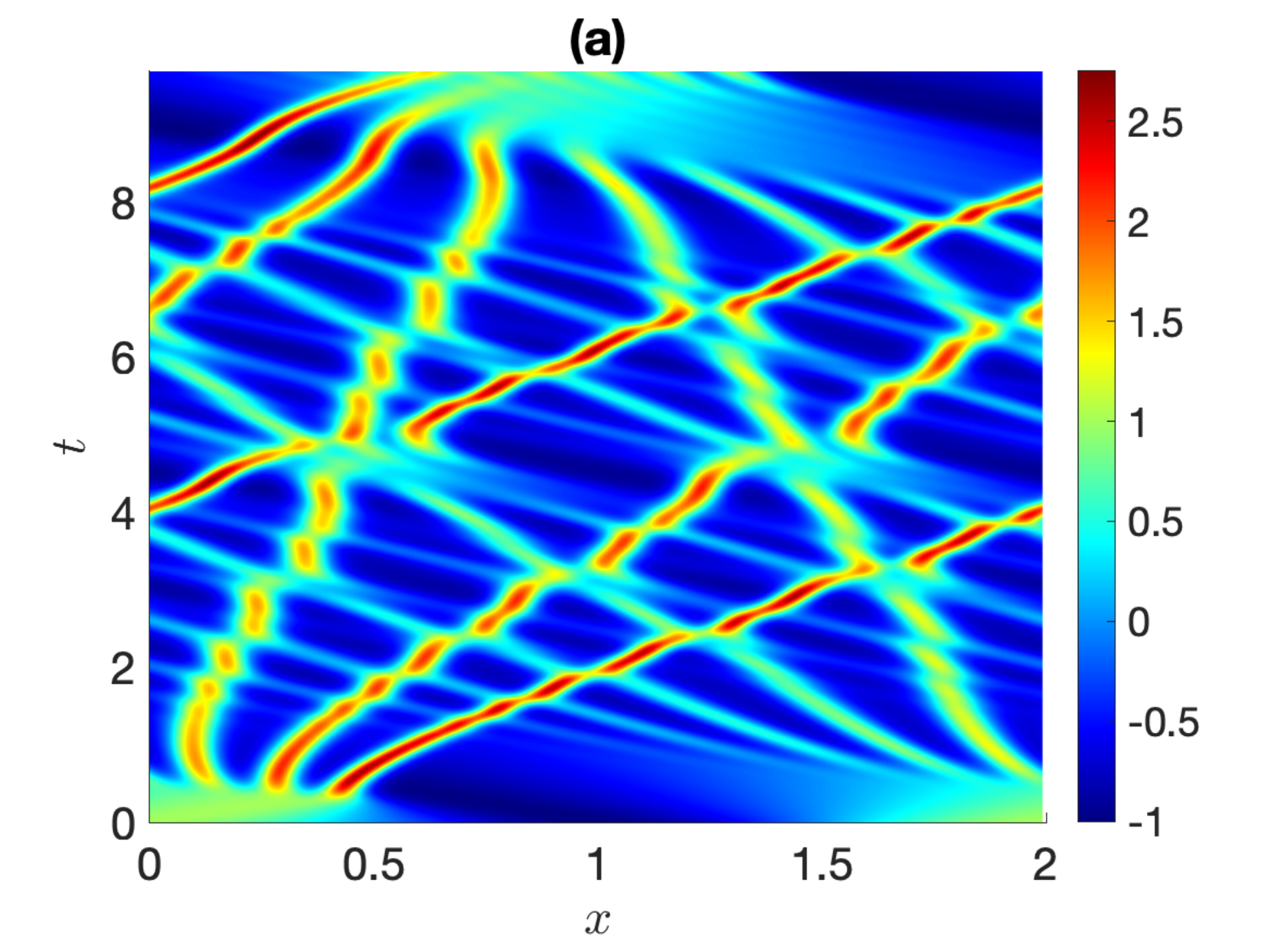}
    \end{subfigure}
    %\hspace{-1.5cm}
    \begin{subfigure}%[h]{0.49\linewidth}
    \centering
    \includegraphics[width=35mm]{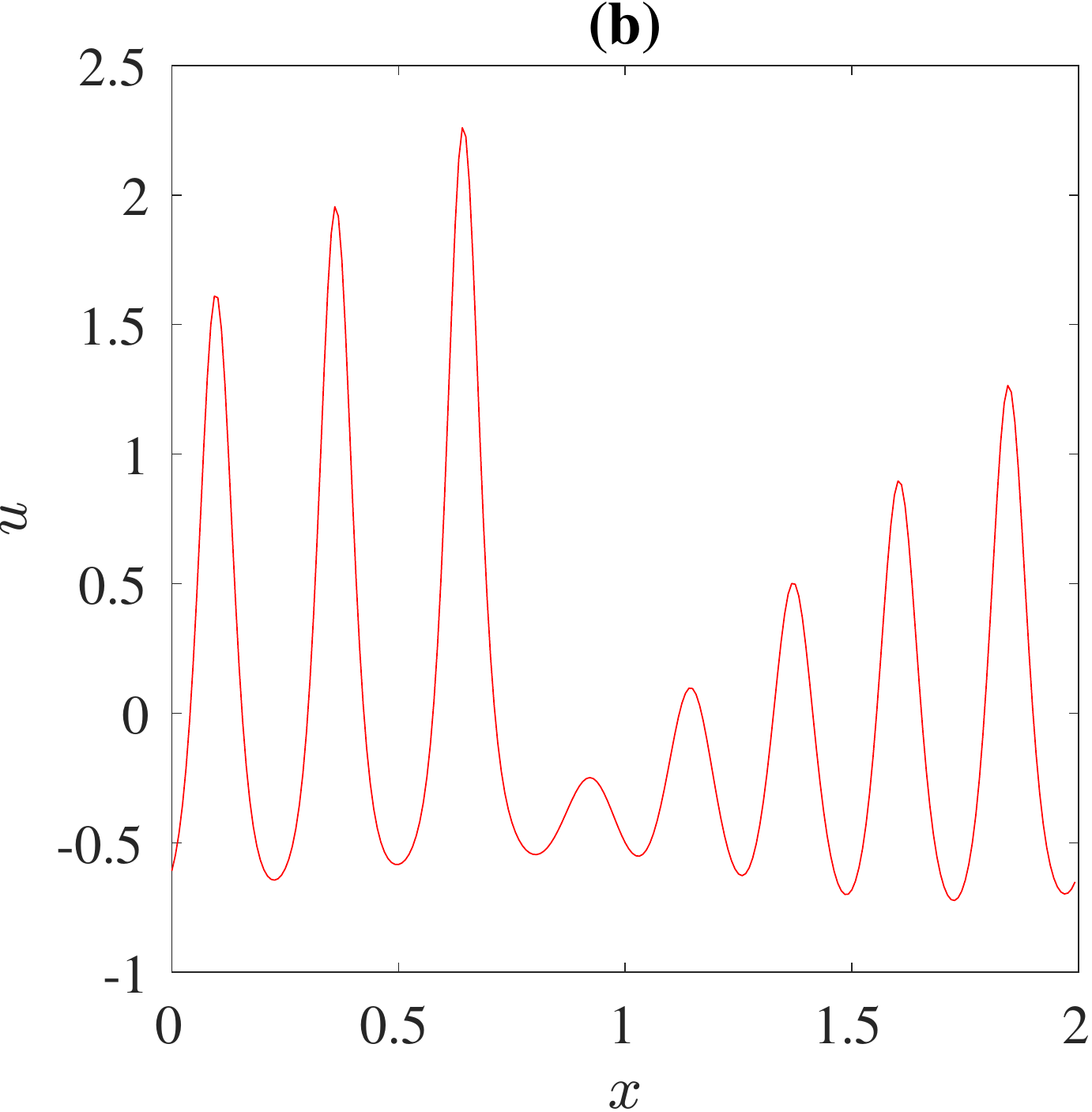}
    \end{subfigure}
    \begin{subfigure}%[h]{0.49\linewidth}
    \centering
    \includegraphics[width=35mm]{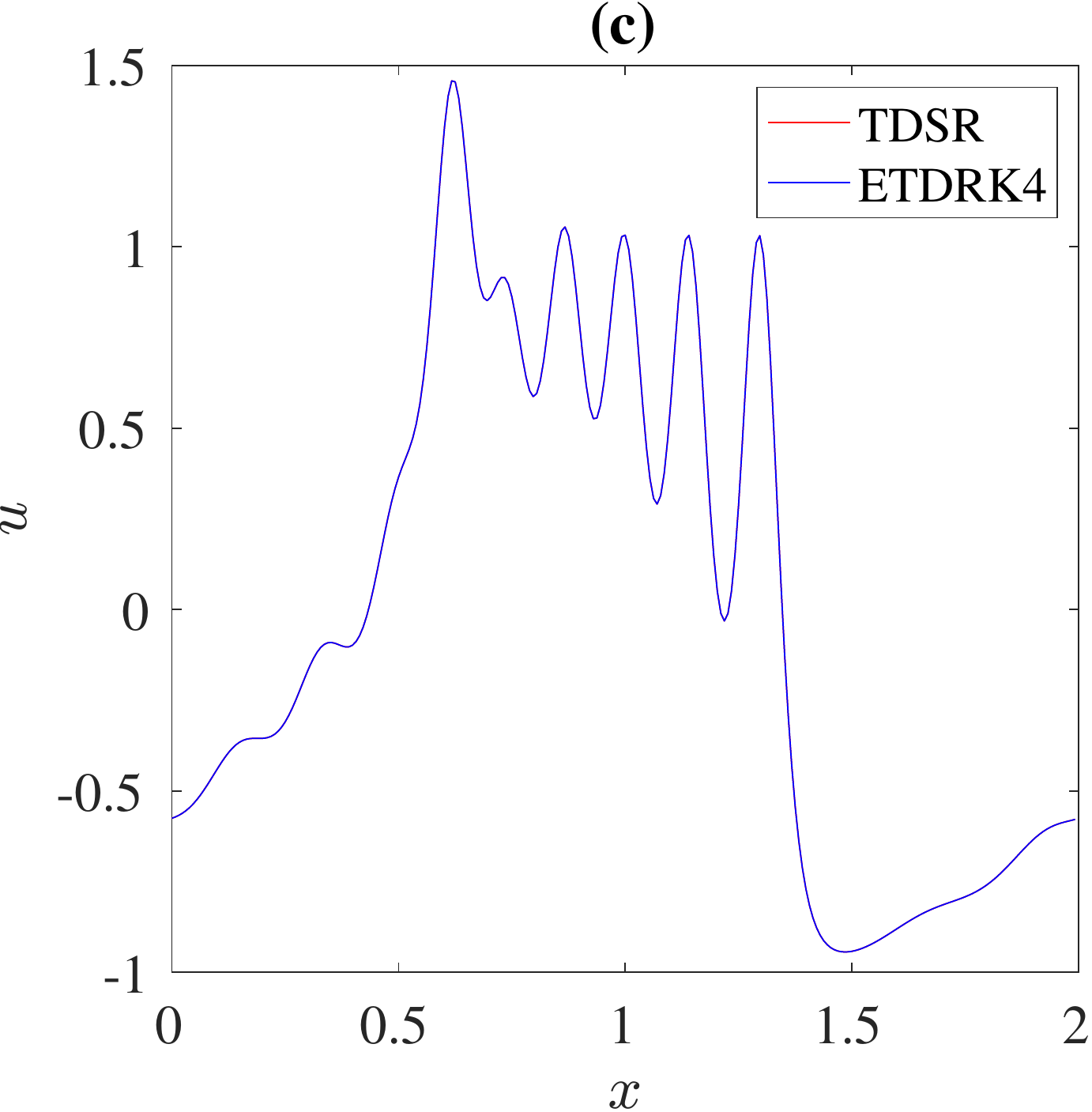}
    \end{subfigure}
    \begin{subfigure}
        \centering
        \includegraphics[width=45mm]{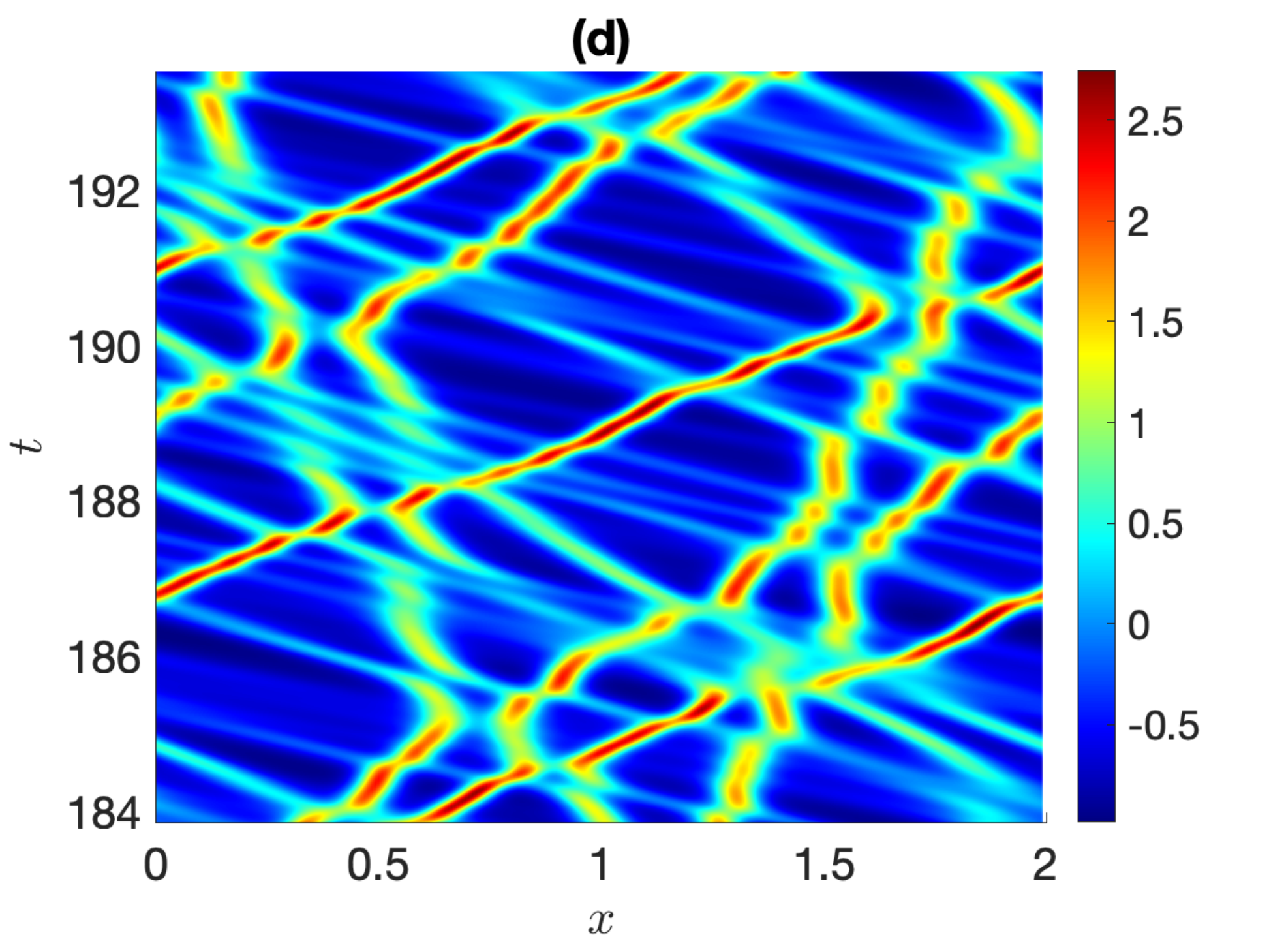}
       % \caption{Caption}
        %\label{fig:my_label}
    \end{subfigure}
    \begin{subfigure}
        \centering
        \includegraphics[width=35mm]{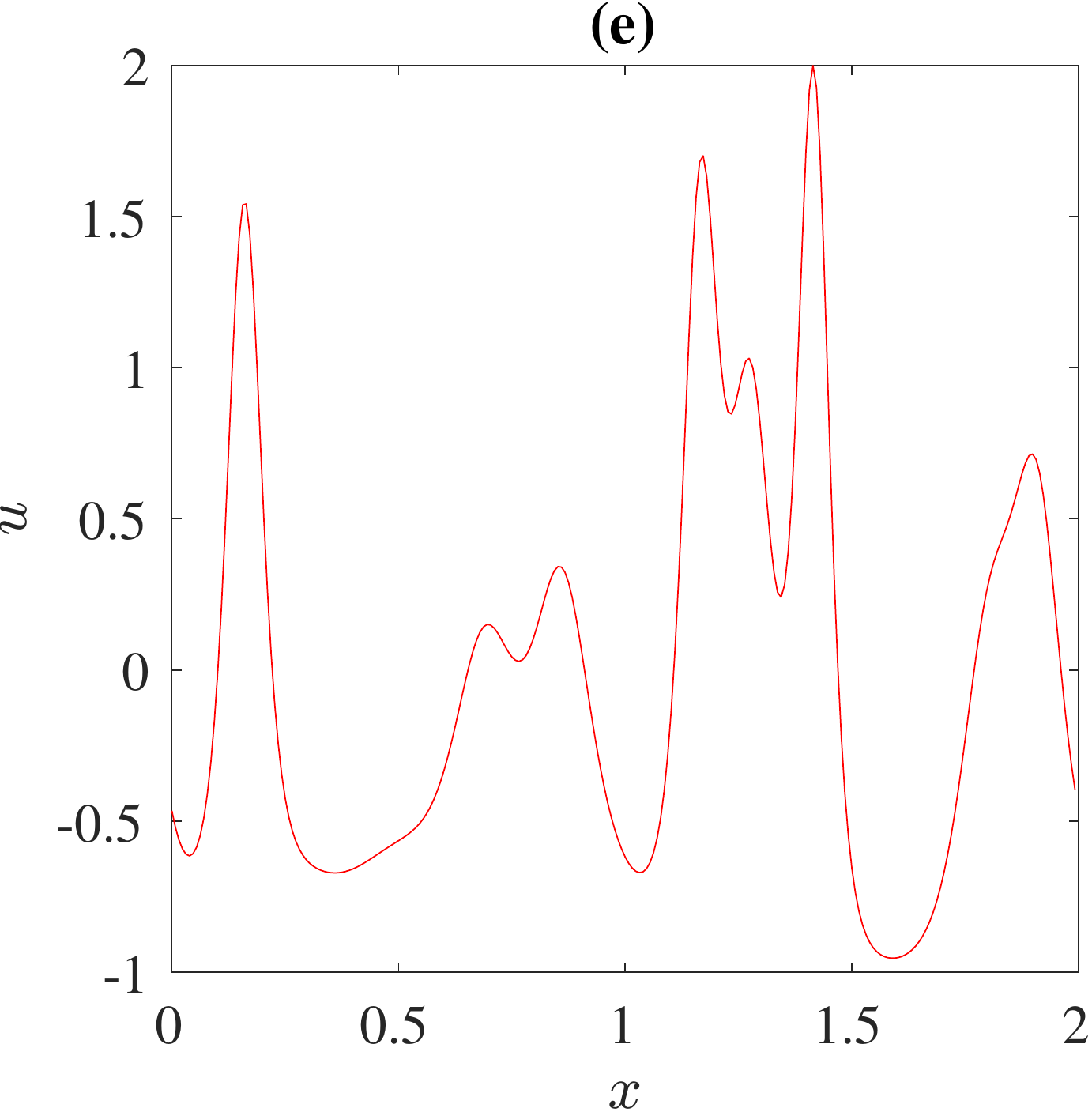}
        %\caption{Caption}
        %\label{fig:my_label}
    \end{subfigure}
   % \end{subfigure}
        \caption{\small (a) Space-time contour plot for the KdV solution with $\alpha =1, \epsilon = 0.022$ and initial condition
    $u(x,0)=\cos(\pi x).$ Other parameters are $\Delta t \approx 0.0008$, $L=2$, $N_S=256$ and the time block size $T_1\approx 0.016$ ({multi-blocking with $M_b$ sub-intervals of equal size}). (b) {The solution at $T=3.6/\pi$ depicting the fission of the initial condition into an eight soliton train.} (c) The solution at one recurrence time $T=t_R=30.4/\pi$ obtained via TDSR Simpson showing good agreement with the fourth-order accurate (in time) ETDRK4 solution. A time step of $\Delta t=0.0001$ (for the ETDRK4 scheme) was used to obtain a solution of comparable accuracy to ours. (d) A space-time contour plot for the solution obtained from the TDSR algorithm using the same spatio-temporal discretization as in (a), over the time span $[19 t_R, 20 t_R]$. The stable numerical simulation produced accurate results, as evidenced by the minor relative errors in the first six conserved quantities (see Fig. \ref{KdV_ZK_LT_D2}). (e) {The solution at $T=20t_R$ obtained using TDSR Simpson.}}
    \label{KdV_ZK_LT_D}
\end{figure}

%%%%%%%%%%%%%%%%%%%%%%%%%%%%%%%%%%%%%%%%%%%%%%%%%%%%
Our goal in this section is to reproduce the well-known numerical results of Zabusky and Kruskal on the KdV equation \cite{zabusky1965interaction} using our algorithm. 
Their simulation displays rich nonlinear dynamics which, as such, represents a challenge for numerical methods as far as the choice of time-steps, long-time accuracy and stability are concerned  \cite{FURIHATA1999181,sanz1982explicit,cui2007numerical,ascher2005symplectic}. 
To this end, we apply the TDSR method on the KdV Eq.~(\ref{kdv-eqn}) with $\alpha =1$ and $\epsilon = 0.022$ 
subject to periodic boundary conditions $u(x+2,t) = u(x,t)$ and initial condition ${u_0(x)=\cos(\pi x)}$. Figure \ref{KdV_ZK_LT_D} (a) shows the Zabusky-Kruskal results while rigorously conserving the momentum $\int_{0}^{2} u^2 dx$.  We compare our findings with those obtained using the ETDRK4 method  \cite{kassam2005fourth}. 
Note that, in order to keep the level of solution accuracy of the ETDRK4 method comparable with that of the TDSR scheme while conserving momentum, the time step $\Delta t_{ETD}$ (of the ETDRK4) has to be about one order of magnitude lesser than its TDSR counterpart ($\Delta t_{TDSR}$). This can be seen by gradually reducing the relative time step ($\Delta t_r=\Delta t_{ETD}/\Delta t_{TDSR}$) while monitoring the relative difference between the two numerical solutions. At $\Delta t_r=0.125$, this difference was seen to drop to $\mathcal{O}(10^{-6})$. %We noted wall clock times of 5 and 4 seconds for the TDSR and ETDRK4, respectively, to arrive at solutions of comparable accuracy.
%%%%%%%%%%%%%%%%%%%%%%%%%%%%%%%%%%%%%%%%%%%%%%%%%%%%%%
Using the same spatio-temporal discretization as before, we checked the TDSR results' fidelity up to 20 recurrence times. This was done, by monitoring the 
time evolution of the first six conserved quantities of the KdV equation with $\alpha=1$ and $\epsilon=0.022$ given by: 
%While the first two conserved quantities correspond to $Q_1=u$ and $Q_2=u^2$, we list the next four below corresponding to $\alpha=1$ and $\epsilon=0.022$ (see \cite{miura1968korteweg} for a list of conserved quantities for the KdV equation):
\begin{align*}
\label{Conserved-quantities-ZK}
% better to multiply by -2!
Q_1 & = u,\qquad Q_2 = u^2\\
Q_3 &= -\frac{u^3}{6}+\frac{\epsilon^2 u_x^2}{2},\qquad 
Q_4 = \frac{u^4}{4}-3\epsilon^2uu_x^2+\frac{9\epsilon^4 u_{xx}^2}{5},\\
Q_5 &= \frac{u^5}{5}-6\epsilon^2u^2u_x^2+\frac{36\epsilon^4uu_{xx}^2}{5}-\frac{108 \epsilon^6u_{xxx}^2}{35},\\
Q_6 &= \frac{u^6}{6}-10\epsilon^2u^3u_x^2+18\epsilon^4u^2u_{xx}^2-5\epsilon^4u_x^4-\frac{108\epsilon^6u u_{xxx}^2}{7}+\frac{120\epsilon^6 u_{xx}^3}{7}+\frac{36\epsilon^8 u_{xxxx}^2}{7}.
\end{align*}

\begin{figure}%[!hb]
%\hspace{0.5cm} 
    \begin{subfigure}%[h]{0.49\linewidth}
    \centering
    %\vspace{-0.9cm}
    %\centering
    %\vspace{-0.9cm}
    %\begin{minipage}{0.6\textwidth}
    \includegraphics[width=120mm]{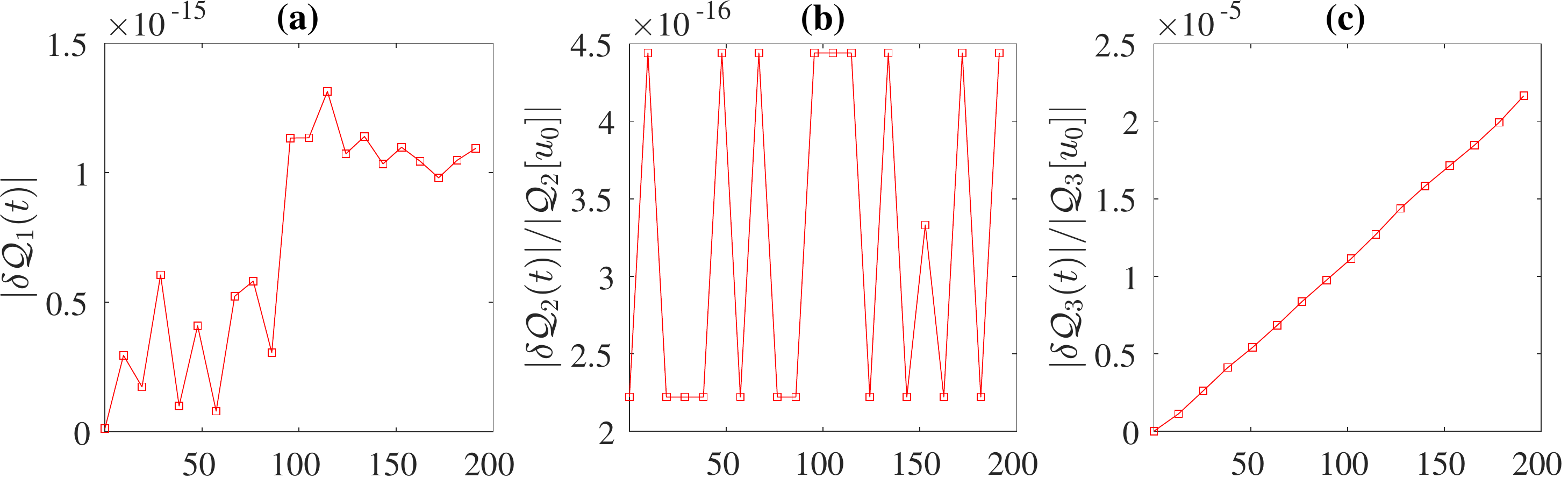}
    \end{subfigure}
    \begin{subfigure}%[h]{0.49\linewidth}
    %\vspace{-0.3cm}
    \centering
    \includegraphics[width=120mm]{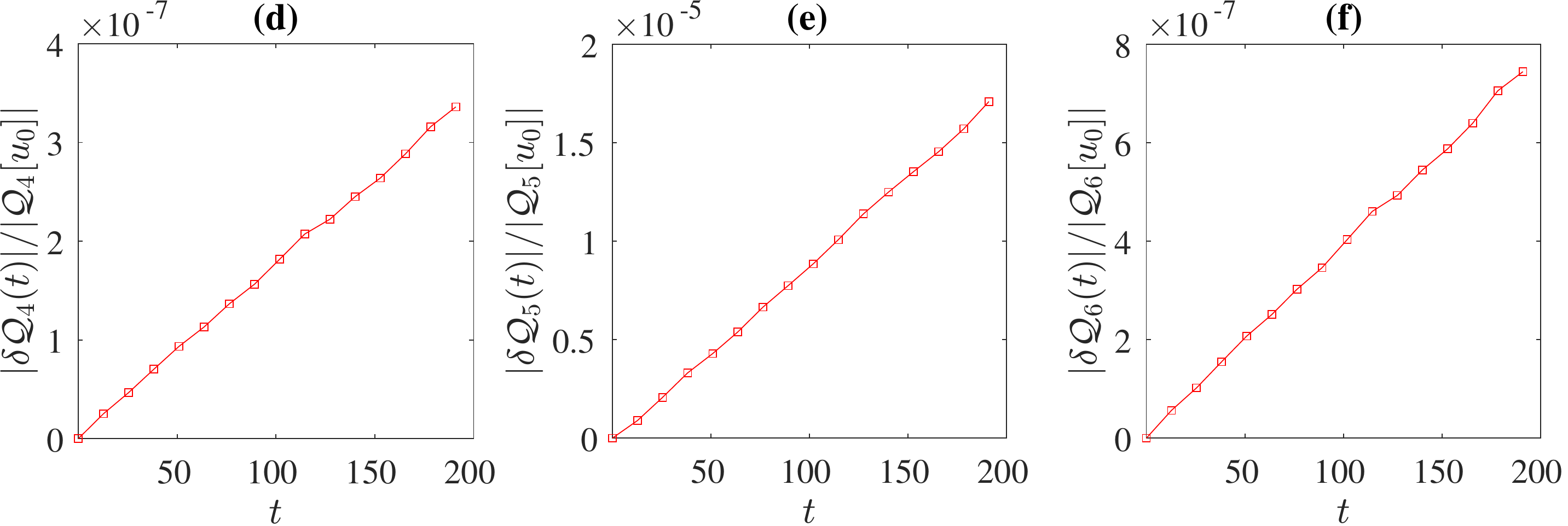}
    \end{subfigure}
   % \vspace{1.0cm}
        \caption{\small {}Errors in conserved quantities for the Zabusky-Kruskal test case monitored over the time-span [0, 20$t_R$]; computational parameters for the simulation can be found in caption of Fig. \ref{KdV_ZK_LT_D}: (a) absolute error in the mass (the initial mass $\mathcal{Q}_1[u_0]=0$), relative errors in (b) momentum ($L^2$-norm), (c) Hamiltonian, (d) fourth conserved quantity ($\mathcal{Q}_4[u_0]$), (e) fifth conserved quantity ($\mathcal{Q}_5[u_0]$) and (f) sixth conserved quantity ($\mathcal{Q}_6[u_0]$). By construction, the relative error in the conservation of momentum is kept near machine precision, while the absolute error in mass remains at $\sim 10^{-15}.$}%The errors in conservation of mass and momentum are found to be close to machine precision for the entire simulation time. }
    \label{KdV_ZK_LT_D2}
\end{figure}

%The simulation took $\sim 100 s$ on our local workstation.
Clearly, the relative error in the momentum remains close to machine precision, while at the same time resulting in the conservation of the mass as well.
{This can be explained by considering the Fourier series representations of the KdV solution $u(x,t)=\sum_{m=-\infty}^{\infty} \hat{u}_m(t) e^{i\pi mx}$, %its square $u^2(x,t)=\sum_{m=-\infty}^{\infty} \hat{U}_m(t) e^{i\pi mx}$ 
and its corresponding auxiliary function $v(x,t)=\sum_{m=-\infty}^{\infty} \hat{v}_m(t) e^{i\pi mx}$. With this at hand, the $(n+1)^{th}$ renormalized Duhamel iterate takes the form:
\begin{equation}
\label{Renorm-Duham}
    \hat{v}_m^{(n+1)}(t)=\frac{1}{R^{(n)}(t)}\bigg(e^{i\epsilon^2 m^3 \pi^3 t}\hat{u}_m(0)-\int_{0}^{t} e^{i\epsilon^2 m^3 \pi^3 (t-\tau)} \frac{i m\pi}{2} \sum_{l=-\infty}^{\infty} \hat{u}_{l}^{(n)}\hat{u}_{m-l}^{(n)} %\hat{U}^{(n)}_m(\tau)  
    d{\tau}\bigg)
\end{equation}
where $(R^{(n)}(t))^2={C_2}/\bigg({\int_{0}^{2}\Big(v^{(n)}(x,t)\Big)^2 dx}\bigg)$. Substituting $m=0$ in Eq.~(\ref{Renorm-Duham}), we obtain:
\begin{equation}
    \hat{v}_0^{(n+1)}(t)=\frac{\hat{u}_0(0)}{R^{(n)}(t)}.
\end{equation}
%This can be explained by first relating $\mathcal{Q}_1[u(x,t)]=2\hat{a}_{0}(t)$ (where $\hat{a}_0(t)$ is the zeroth Fourier modal coefficient in the Fourier series representation $u(x,t)=\sum_{n=-\infty}^{\infty} \hat{a}_n(t) e^{i\pi kx}$), and thereafter viewing the renormalized Duhamel iteration in Fourier space for $n=0$:

Using the identity $\hat{u}_0^{(n+1)}(t)= R^{(n+1)}(t)\hat{v}_0^{(n+1)}(t)$ we find
\begin{equation}
    \hat{u}_0^{(n+1)}(t)=\sqrt{\frac{{\int_{0}^{2}\Big(v^{(n)}(x,t)\Big)^2 dx}}{{\int_{0}^{2}\Big(v^{(n+1)}(x,t)\Big)^2 dx}}}\hat{u}_0(0).
\end{equation}
For the initial condition considered here $u_0(x)=\cos(\pi x)$, one finds $\hat{u}_0(0)=0$, resulting in $\hat{u}_0^{(n+1)}(t)\equiv 0$, i.e., the mass is also preserved at \textit{every} Duhamel iterate.
Additionally, the relative errors of the other four conserved quantities are within $\mathcal{O}(10^{-5})$. 
In particular, some aspects of our scheme (such as solution accuracy) outperforms other well known conservative numerical methods to simulate the KdV equation, such as the one developed in \cite{cui2007numerical}. In  \cite{cui2007numerical}, the authors develop an operator splitting scheme in conjunction with a finite volume spatial discretization to locally conserve the mass and momentum. 
While their scheme demonstrates impressive long time stability properties (Zabusky-Kruskal dynamics), there were some phase errors at $T=20t_R$, that arise from the global (absolute) errors in the conservation of the Hamiltonian ($\sim 10^{-3}$).
%Thus, the drawback of their idea is the lack of generalization to an arbitrary number of conservation laws.

%Within each finite volume, a linear (spatial) interpolant was used at the \textit{reconstruction} step. Here, the authors proposed to compute the slope of the interpolant by enforcing  the momentum average of the reconstructed solution to be equal to the numerical momentum within the finite volume.  Crucially, they showed that the slope determination is \textit{feasible}, thus allowing for the extension to incorporate both the local mass and momentum conservation laws of the KdV. 
 Presently, it is unclear if there are other finite volume based schemes capable of incorporating more than two conserved quantities for the KdV equation.} {Also, it is worth to mention that the ETDRK4 scheme (which does not conserve momentum at the local or global level for the KdV), suffers from numerical instabilities (for $\Delta t\approx 0.0008$) and for time integration up to $T=20t_R$.}
%Finally, we also mention that we had to halve the time-step for the ETDRK4 scheme to avoid numerical instability.} \\

{Finally, we examine the performance of the TDSR method by measuring the local errors in the mass and momentum. Specifically, we compute the errors for mass and momentum at end time $T=20t_R$ respectively defined by
\begin{equation}
\label{local-mass-err}
    \mathcal{E}_1(x,T)=u_t+\alpha uu_x+\epsilon^2 u_{xxx},
\end{equation}
\begin{equation}
\label{local-mom-err}
    \mathcal{E}_2(x,T)={(u^2/2)}_t+\bigg(\frac{\alpha}{3}u^3+\epsilon^2\Big(uu_{xx}-{u_x^2}/2\Big)\bigg)_{x}.
\end{equation}
{The time derivatives are computed using the fourth-order backward differentiation formula \cite{quarteroni2010numerical}:
\begin{align}
    u_t(x_m,T)\approx \frac{\Big(25 u(x_m,T)-48 u(x_m,T-\Delta t)+36 u(x_m,T-2\Delta t)-16u(x_m,T-3\Delta t)+3u(x_m,T-4\Delta t)\Big)}{12\Delta t},
    %\mathcal{E}_1(x_m,T)&\approx u(x_{m},T)-\frac{48}{25}u(x_m,T-\Delta t)+\frac{36}{25} u(x_m,T-2\Delta t)-\frac{16}{25} u(x_m,T-3\Delta t)+\frac{3}{25} u(x_m,T-4\Delta t)\\\nonumber&-\left.\frac{12}{25}\Delta t \Big(\mathcal{L}(u)+\mathcal{N}(u)\Big)\right\vert_{x_m,T}.
\end{align}
whereas, all the spatial derivatives were computed to spectral accuracy with the use of fast Fourier transforms.
}
%For $t\in [T_{i-1},T_i]$ (the $i^{th}$ contractive time sub-interval), the local error in the mass conservation law is computed by recasting Eq.~\ref{local-mass-err} in an integral form:
%\begin{equation}
  %  \label{local-mass-err-integ}
  % \mathcal{L}\mathcal{E}_1(x,t)=u(x,t)-e^{-(t-T_{i-1})\epsilon^2(\partial^3/\partial x^3)}[u(x,T_{i-1})]-\int_{T_{i-1}}^{t} e^{-(t-\tau)\epsilon^2(\partial^3/\partial x^3)}\Big(-\alpha uu_x\Big) d{\tau}
%\end{equation}
%and the local error in momentum conservation law is computed by recasting Eq.~\ref{local-mom-err} in integral form as well:
%\begin{equation}
 %   \label{local-mom-err-integ}
 %   \mathcal{L}\mathcal{E}_2(x,t)=u^2(x,t)-u^2(x,T_{i-1})-\int_{T_{i-1}}^{t} -\bigg(\frac{2\alpha}{3}u^3+\epsilon^2\Big(2uu_{xx}-{u_x^2}\Big)\bigg)_x d{\tau}
%\end{equation}
\begin{figure}
    \centering
    \includegraphics[width=0.58\linewidth]{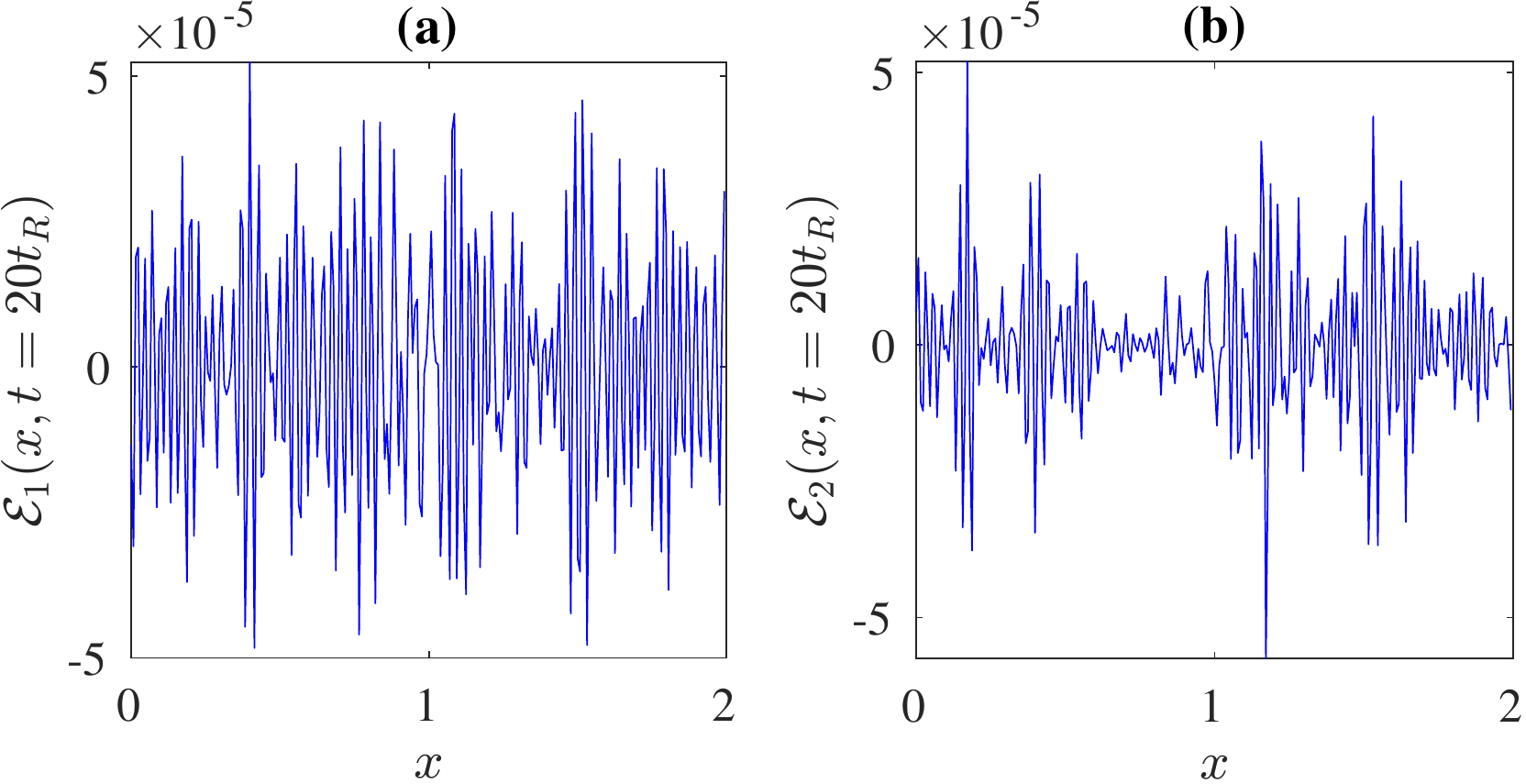}
    \caption{{A snapshot of the local errors in (a) conservation of mass ($\mathcal{E}_1(x,t)$) and (b) conservation of momentum ($\mathcal{E}_2(x,t)$) as a function of $x$ at time $T=20t_R$.}}
    \label{Local-mass-mom-cons}
\end{figure}
 Fig.~\ref{Local-mass-mom-cons} shows the variations of $\mathcal{E}_1$ and $\mathcal{E}_2$ as a function of $x$ at time $T=20t_R$. Remarkably, the local errors remain relatively small even over such long time intervals.
}
%%%%%%%%%%%%%%%%%%%%%%%%%%%%%%%%%%%%%%%%%%%%%%%%%%%%%%%%

%\bigskip
%%%%%%%%%%%%%%%%%%%%%%%%%%%%%%%%%%%%%%%%%%%%%%%%%%%%%%%
\subsection{Travelling waves for the Allen-Cahn equation.} 
%%%%%%%%%%%%%%%%%%%%%%%%%%%%%%%%%%%%%%%%%%%%%%%%%%%%%%%
In this section, we apply the TDSR method on the Allen-Cahn equation, a prototypical reaction-diffusion type equation that arises in material science \cite{ALLEN19791085}. It is given by
\begin{equation}
\label{AC-neumann}
 u_t = D u_{xx} + \gamma (u-u^3),
 \end{equation}
where $D>0$ is the diffusion coefficient and $\gamma$ measures the strength of reaction. The Allen-Cahn equation is dissipative in nature. In fact when subject to homogeneous Dirichlet/Neumann boundary conditions (considered in this paper), multiplying Eq.~(\ref{AC-neumann}) by $2u$ and integrating over the whole domain leads to
the dissipation rate Eq.~(\ref{n-diss-rates}) with density and flux:
\begin{equation}
\label{AC-neumann-diss}
 \rho = u^2,\;\;\;\;\; F= 2Du^2_x - 2 \gamma  (u^2-u^4).
 \end{equation}
In this case, as we shall see later, the renormalization factor $R(t)$ obeys a nonlinear ordinary differential equation.
We present numerical results on two canonical problems associated with the Allen-Cahn equation: (i) dynamics of traveling waves and (ii) observation of meta-stable dynamics. Both examples represent a departure from the periodic case for which the linear 
operator $\mathcal{L}$ is diagonalizable. Indeed, for the Allen-Cahn equation, the discrete representation 
of $\mathcal{L}$ is now dense as is the case when using spectral differentiation matrices. Thus, the semi-group $\exp (t{\mathcal{L}})$ forms a 
rank-3 tensor.
We implement the TDSR method on the Allen-Cahn Eq.~(\ref{AC-neumann}) to compute traveling waves with diffusion coefficient $D=1$ and large reaction parameter $\gamma$ that scales like $\epsilon^{-2}$, with $\epsilon\ll 1.$
 Equation \eqref{AC-neumann} is subject to the initial condition $u(x,0) = 0.5-0.5\tanh\left(x/(2\sqrt{2}\epsilon)\right)$ and Neumann boundary conditions: 
$u_x (x,t) \rightarrow  0$ as $|x|\rightarrow \infty.$ Interestingly enough, the Allen-Cahn Eq.~(\ref{AC-neumann}) admits an exact travelling wave solution given by $u(x,t)=0.5-0.5\tanh \xi /(2\sqrt{2}\epsilon),\;\; \xi = x - 3t/(\sqrt{2}\epsilon)$ \cite{jeong2016comparison}. It is the aim of this section to reproduce this exact solution using the TDSR while enforcing the dissipation rate equation given in (\ref{n-diss-rates}) and (\ref{AC-neumann-diss}). We proceed by
substituting the ansatz $u(x,t)=R(t) v(x,t)$ into Eq.~(\ref{AC-neumann}); multiply by $2u$ and integrate the resulting system over the whole computational domain to obtain a first order dynamical system for the variable $p(t) \equiv  r(t) R^2(t)$: 
%%%%%%%%%%%%%%%%%%%%%%%%
\begin{equation}
    \label{ODE-AC-Renorm}
    \frac{dp}{dt} = (-a(t)+2\gamma) p- b(t) p^2,
\end{equation}
%%%%%%%%%%%%%%%%%%%%%%%%%
where the expressions for the time-dependent coefficients $r(t)$, $a(t)$ and $b(t) >0$
%Appendix \ref{AC-renorm-appendix}.
are given by (here $\Omega$ denotes the spatial domain of the Allen-Cahn equation)
\begin{equation}
\label{Allen-cahn-dissipation}
\begin{aligned}
r(t) = \int_{\Omega} v^2(x,t) dx, \quad 
a(t) = \frac{2D\int_{\Omega} v_x^2(x,t) dx}{r(t)}, \quad
b(t) =  \frac{2\gamma\int_{\Omega} v^4(x,t) dx}{r^2(t)}.
\end{aligned}
\end{equation}
%In fact, these time-dependent variables are functionals that depend on the auxiliary wave function $v(x,t).$

%%%%%%%%%%%%%%%%%%%%%%%%%%%%%%%%%%%%%%%%%%%%%%%%%
\begin{figure}
\begin{subfigure}%{.49\textwidth}
\centering
\includegraphics[width=60mm]{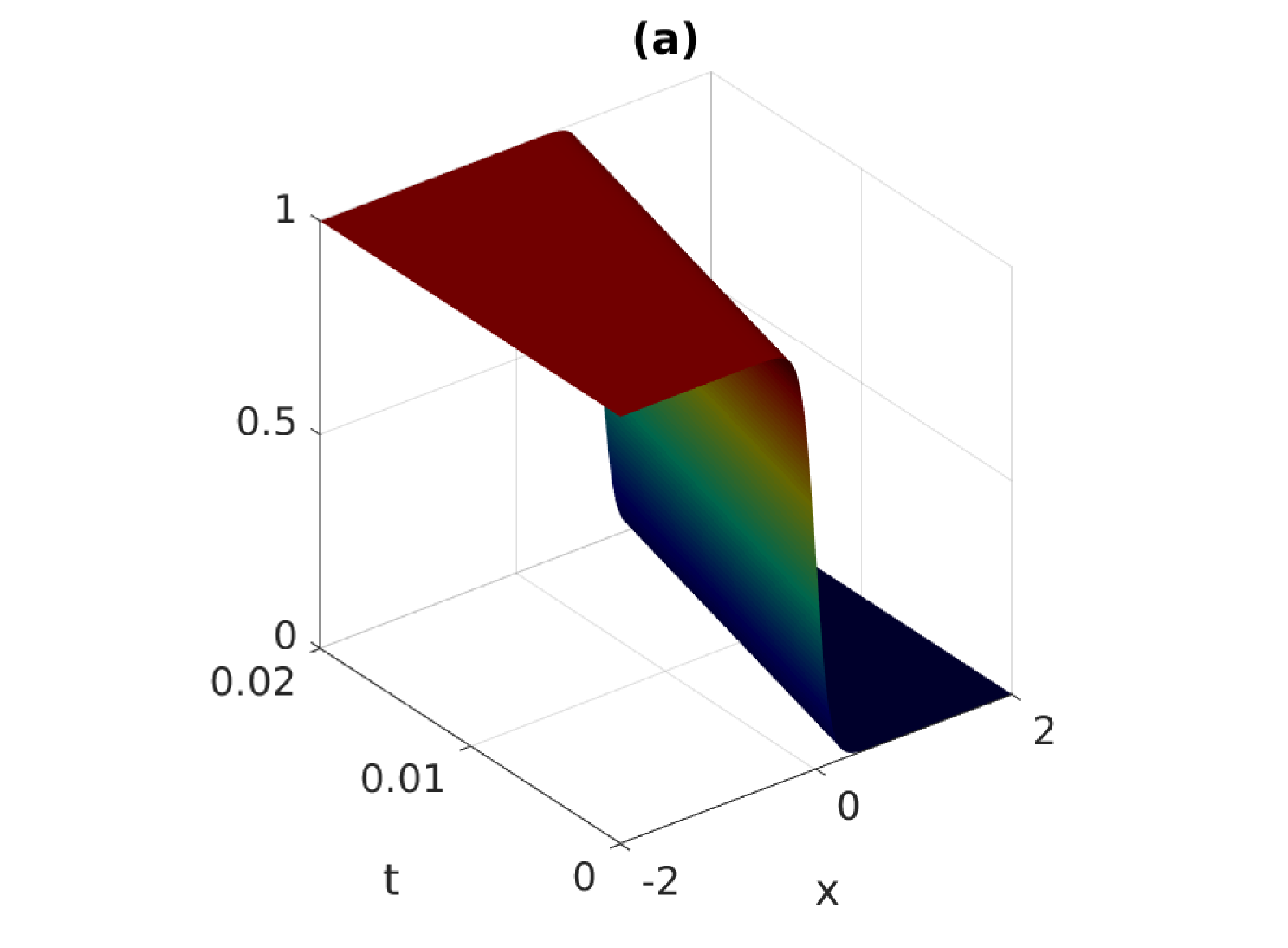}
\end{subfigure}
\begin{subfigure}%{.49\textwidth}
\centering
\includegraphics[width=40mm]{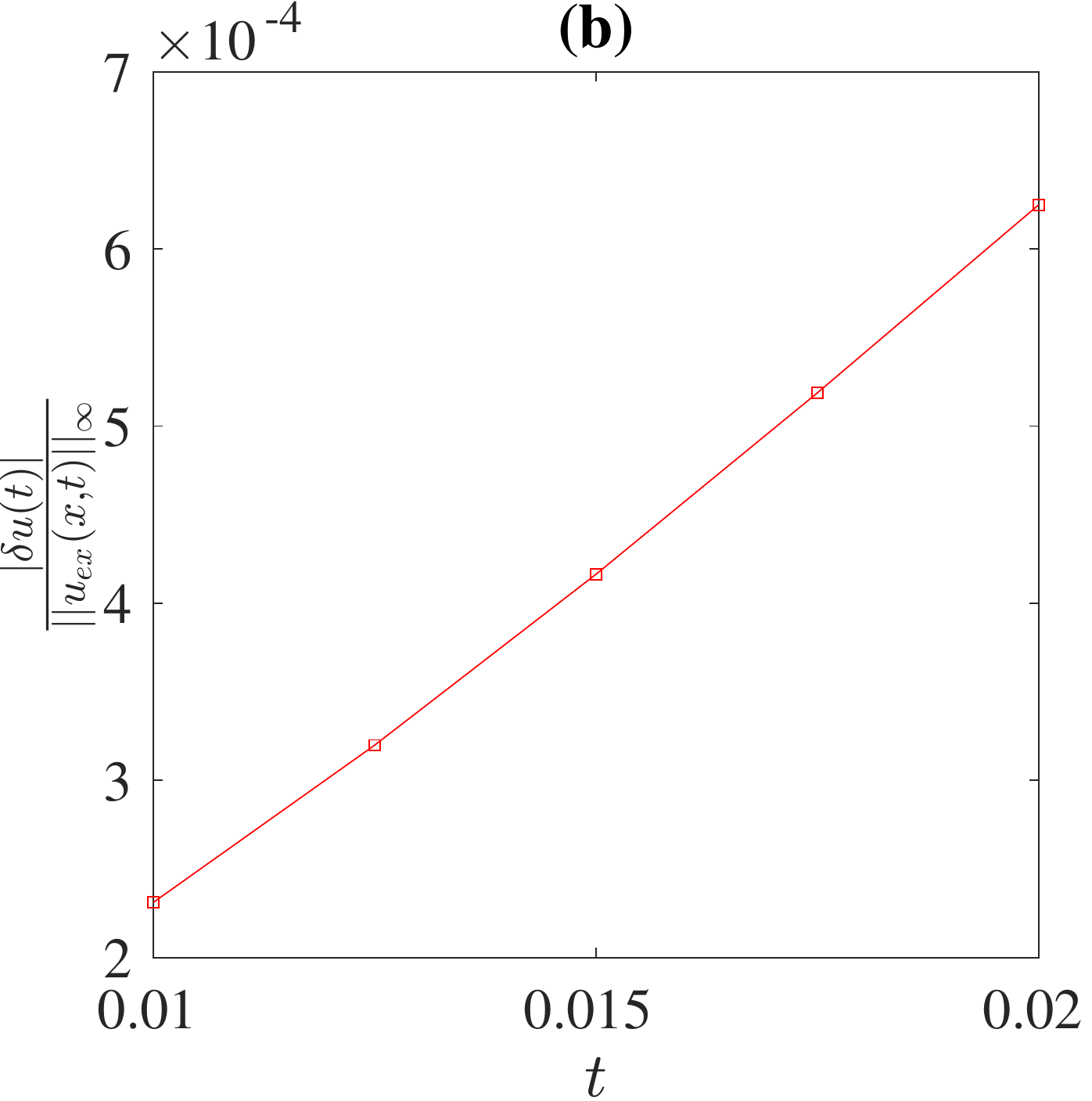}
\end{subfigure}
%\hfill
\caption{\small(a) Travelling wave solution for the Allen-Cahn equation with simulation parameters: $L=4$, $\epsilon=0.05$, $N_S=1024$, $\Delta t= 0.000125$, and end time $T=0.02$. Solution advanced to the final time with the idea of multi-blocking with sub-interval size $T_1=0.01$. The error in solution was restricted to $6.3\times 10^{-4}$ with this parameter choice. (b) The absolute error, in time, between the TDSR numerical and the exact solutions in the time interval $[0.01, 0.02]$.} 
\label{Travelling_wave_results}
\end{figure}
%%%%%%%%%%%%%%%%%%%%%%%%%%%%%%%%%%%%%%%%%%%%%%%%%%

The presence of the large coefficient $\gamma\sim \epsilon^{-2}$ in Eq.(\ref{ODE-AC-Renorm}) causes the differential equation to become stiff, thus severely limiting the choice of time-steps. With this in mind, we use an implicit scheme (such as Crank-Nicolson) to time-step (\ref{ODE-AC-Renorm}) \cite{quarteroni2010numerical}. The coefficients $r(t)$, $a(t)$ and $b(t)$ are computed to spectral accuracy using Clenshaw-Curtis quadrature method \cite{trefethen2000spectral}. 
 %%%%%%%%%%%%%%%%%%%%%%%%%%%%%%%%%%%%%%%%%%%%%%%%%%%%%%%%%%%%
 The spatial domain $\Omega$ is truncated, in which case, $x\in [-L/2,L/2]$ and taking advantage of the spatial decay of $u_x$ to enforce homogeneous Neumann boundary conditions. Several remarks are in order: (i) The Chebyshev (Chebyshev-Lobatto) series representation is originally developed for functions defined on the interval $ [-1,1].$ It can, nonetheless be applied on the interval $[-L/2,L/2]$ using a linear transformation. (ii) We incorporate the homogeneous Neumann boundary conditions following similar procedure as outlined in \cite{liu2014exponential}. For example, the second derivative is represented by the matrix product $D D_0$ with $D$ being the standard first order spectral differentiation matrix while $D_0$ is the first order differentiation matrix whose first and last rows have been replaced by the zero-vector respectively. This implicitly enforces the underlying boundary conditions. The numerical results obtained in this section are summarized in Fig. \ref{Travelling_wave_results} where a surface plot describing the time evolution of the traveling wave as well as its numerical accuracy are shown.  %(b), where we track the time evolution of the absolute error in the TDSR numerical solution.

\bigskip
%%%%%%%%%%%%%%%%%%%%%%%%%%%%%%%%%%%%%%%%%%%%%%%%%%%%%%%%%%
\subsection{Meta-stable dynamics of the Allen-Cahn equation.}
\begin{figure}
\begin{subfigure}%{.49\textwidth}
\centering
\includegraphics[width=50mm]{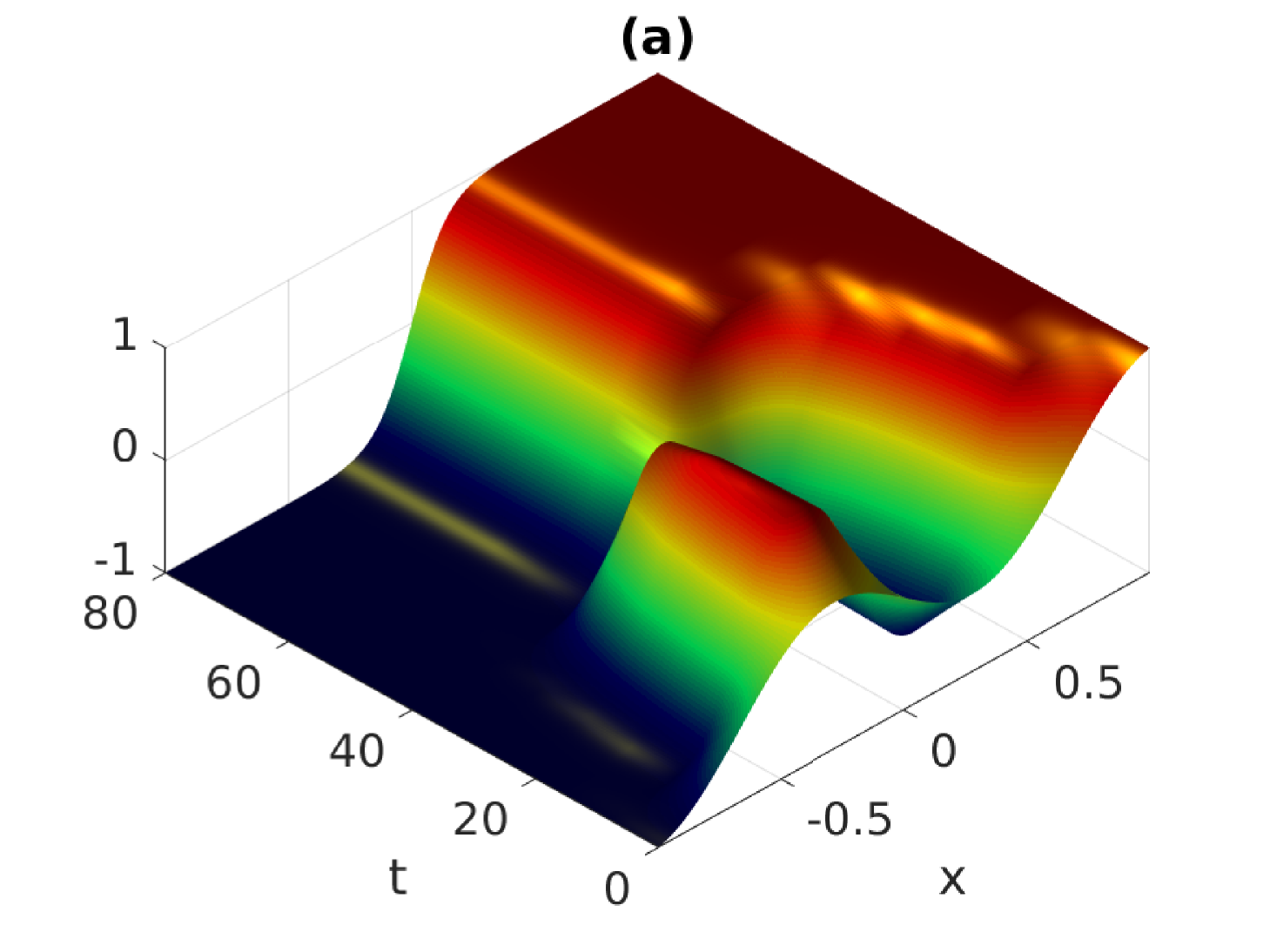}
\end{subfigure}
\begin{subfigure}%{.49\textwidth}
\centering
\includegraphics[width=50mm]{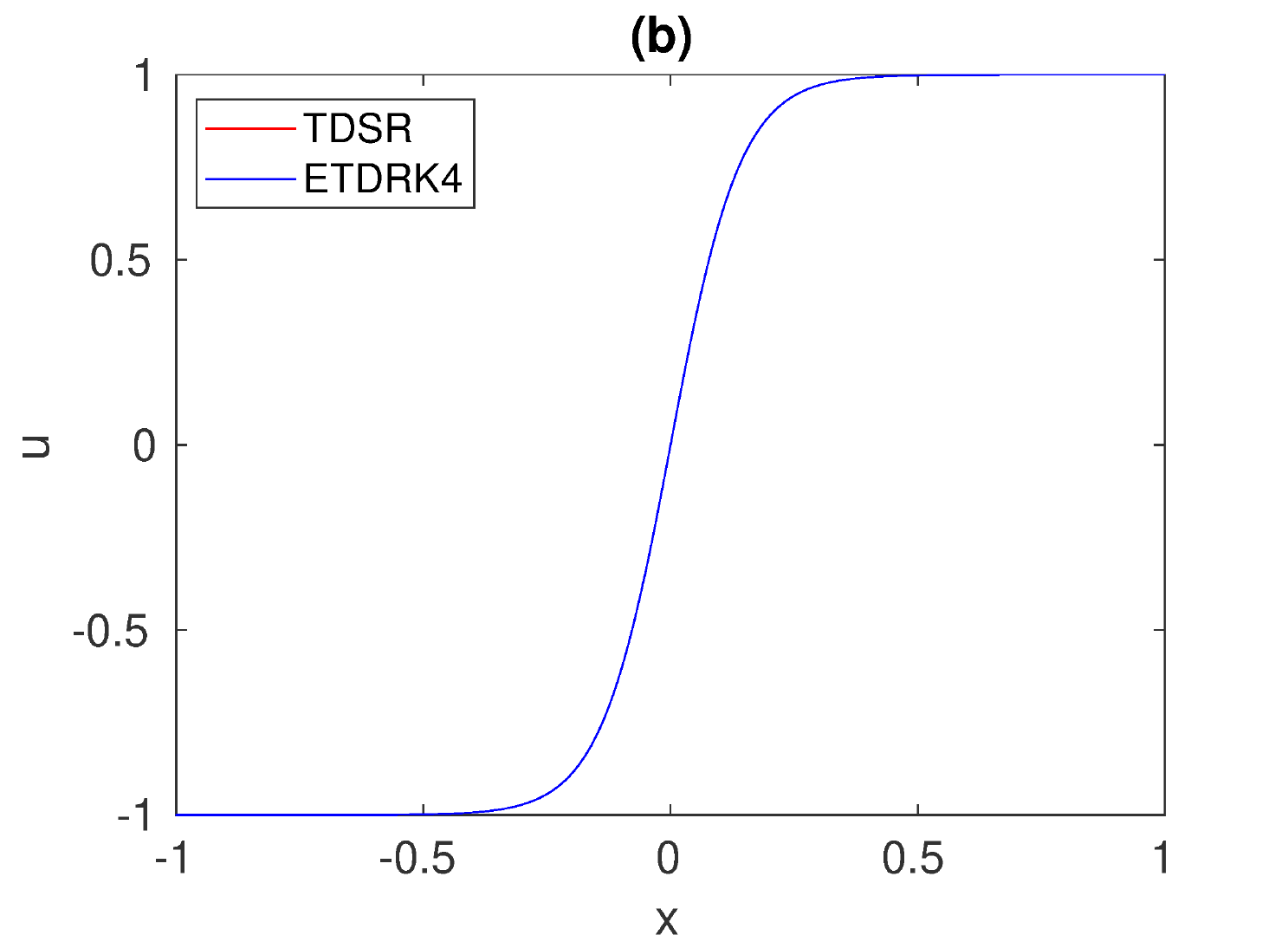}
\end{subfigure}
%\hfill
\caption{\small (a) Spatio-temporal field distribution illustrating the meta-stable dynamics. Simulations were performed using the idea of multi-blocking with block size $T_1=8$ and final time $T=80$. Computational parameters:
$N_S=256$, $\Delta t \approx 0.016$. (b) Comparison with the ETDRK4 scheme (at end time), showing a good agreement with the TDSR trapezoidal result. }%The TDSR trapezoidal numerical solution showing good agreement with the results from } 
\label{MetaStable_results}
\end{figure}

%%%%%%%%%%%%%%%%%%%%%%%%%%%%%%%%%%%%%%%%%%%%%%%%%%%%%%%%%%%%%%%%%%%
%\begin{figure}
 %   \centering
 %   \includegraphics[width=0.5\linewidth]{Fig14.pdf}
 %   \caption{\redc{Figure illustrating the role of the renormalization factor (plotted over the $6^{th}$ time block) in enforcing convergence. Without renormalization, each subinterval size had to be restricted to $T_1=5$ in order to ensure the convergence of the Duhamel iterations.}}
 %   \label{fig:Renormalization-allen-cahn}
%\end{figure}
%%%%%%%%%%%%%%%%%%%%%%%%%%%%%%%%%%%%%%%%%%%%%%%%%%%%%%%%%%%%%%%%%%%
Our last example is concerned with the dynamics of a meta-stable state associated with the Allen-Cahn 
Eq.~(\ref{AC-neumann}) with parameters $D =0.01, \gamma=1$ and subject to the boundary 
conditions $u(-1,t)=-1,~ u(1,t)=1$ and initial condition $u(x,0) = 0.53 x+0.47 \sin(-1.5 \pi x).$
This test bed case is particularly interesting since the dynamics of an initial hump is observed to be meta-stable, i.e., it remains unchanged over long time, before abruptly vanishing. This type of rapid change in the wave profile over short time scales inevitably creates numerical challenges. Such dynamic metastability was numerically observed by Kassam and Trefethen \cite{kassam2005fourth} using the modified ETDRK4 scheme.  %example was first considered with the modified ETDRK4 scheme by x \cite{kassam2005fourth}. %\Revone{The }% a sudden and rather abrupt change in the solution at around $t\approx 45$ was noted in the original simulation.
In this section, we demonstrate the robustness of our TDSR method by reproducing this type of abrupt transition from a meta-stable state to another stable wavefunction profile.
%Thus, we envision a prominent role of the renormalization factor to ensure the convergence of our fixed point iterations. 
Here, the renormalization factor $R(t)$ is governed by the same nonlinear ordinary differential equation Eq.~(\ref{ODE-AC-Renorm}) with the exception that now  the parameters are $D =0.01, \gamma=1$ with a spatial domain $[-1,1]$. 
Few remarks are in order: (i) To simplify the computation, we first homogenize the boundary conditions  by deriving a new evolution equation on which the TDSR method is implemented. (ii) To impose Dirichlet boundary conditions, the operator $d^2/dx^2$ is approximated by $D^2$ where $D$ is the first order spectral differentiation matrix \cite{trefethen2000spectral}. 
The time evolution of the Allen-Cahn front is shown in 
Fig.~\ref{MetaStable_results} (a). As expected, our method is indeed capable of reproducing those well known results. We have also compared our results with those  obtained using the ETDRK4 scheme and found good agreement (see Fig.~\ref{MetaStable_results} (b) for a comparison at end time). 

\bigskip 
\bigskip
%%%%%%%%%%%%%%%%%%%%%%%%%%%%%%%%%%%%%%%%%%%%%%%%%%%%%%%%%%%%%%%%%%%%%%%%%%%%%
\section{Numerical Implementation of TDSR with multi-conservation laws}
\label{Multiple-conservation-laws-section}
%%%%%%%%%%%%%%%%%%%%%%%%%%%%%%%%%%%%%%%%%%%%%%%%%%%%%%%%%%%
So far we have addressed several cases where a single conservation law is ``injected" into the numerical simulations. In this section,
we shall present results when multiple conservation laws are enforced. There are three choices that we considered: conservation of (i) mass and momentum;  (ii) mass and Hamiltonian ; (iii) mass, momentum and Hamiltonian. All numerical results reported in this section are for the KdV equation \eqref{kdv-eqn}, subject to rapidly decaying boundary conditions with $\alpha=6$ and $\epsilon =1$. 

{We remark that the pseudo initial conditions $f_j({\bf x})$, $j=1,2,\cdots,N$, introduced in the TDSR formulation (see Eqs. \ref{eqn-phi-1}-\ref{eqn-phi-j}) are crucial for the success of the method. They greatly control the scheme's convergence and allow the renormalization factors to act as ``expansion" coefficients. Our numerical tests strongly indicate that taking none of the pseudo-initial conditions be identically equal to zero, (albeit satisfying the underlying boundary conditions), in order for the scheme to converge. With this in mind, a natural and important issue that immediately arises is how to choose them? 
Our extensive numerical experiments seem to suggest that the natural choice $f_j({\bf x}) \equiv \alpha_j u_0(\textbf{x})$ with non-zero constant $\alpha_j$'s does not lead to convergence. However, several other possible choices for $f_j$ are given by $f_j({\bf x}) = \alpha_j(\textbf{x}) u_0(\textbf{x})$ where $\alpha_j(\textbf{x})$ are spatially localized functions.  
For evolution equations in $(1+1)D$ subject to rapidly decaying boundary conditions (such as the KdV, NLS, mKdV), the algorithm seems to converge (at least over some time interval) when the first $(N-1)$ $f_j({x})$ are chosen to belong to the class of bell-shaped, sign-definite functions. For example, $f_j\in \{{\rm sech}({ x}), e^{-x^2}, {\rm sech}^2({x})\}$ for $j=1,2,\cdots (N-1)$, while the entire set of pseudo-initial conditions is required to satisfy the normalization condition in Eq.~\eqref{fj-condition}. A full characterization on the choice of the pseudo-initial conditions ${f_j({\bf x})}$ is the subject for future work.}

\bigskip
%%%%%%%%%%%%%%%%%%%%%%%%%%%%%%%%%%%%%%%%%%%%%%%%%%%%%%%
\subsection{Conservation of mass and momentum}
\begin{figure}
    \centering
    \includegraphics[width=150mm]{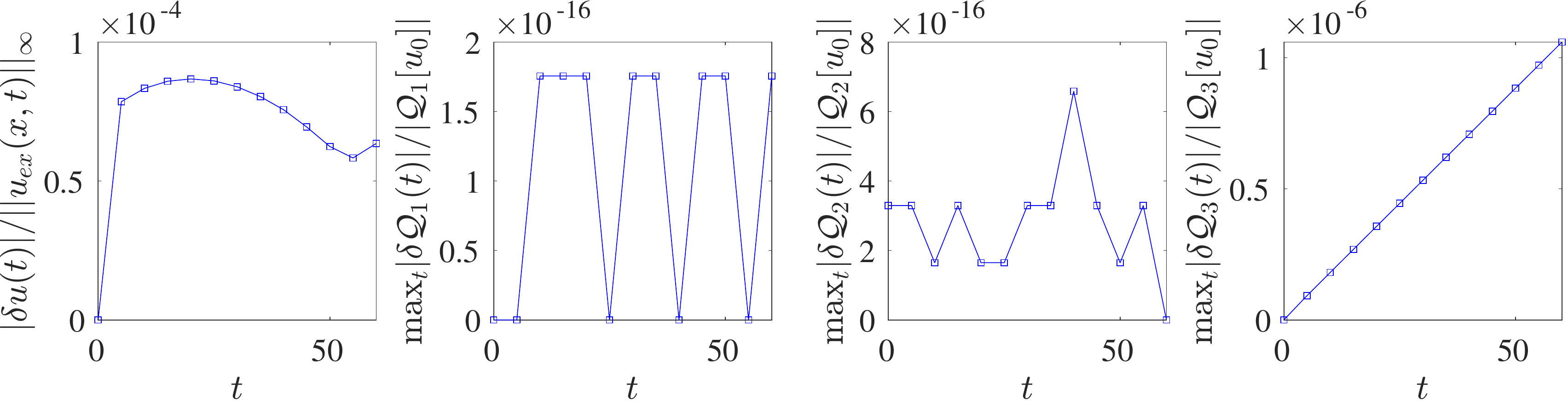}
    \caption{\small (a) Time evolution of the {relative} error between the numerically obtained solution for the KdV equation compared to its exact solution given in Eq.~(\ref{solution-accuracy}). Time evolution of the conserved quantities given in
Eq.~(\ref{cons-law-convergence}): mass (b), momentum (c) and Hamiltonian (d). Note that the relative errors in conservation of mass and momentum stay close to machine precision. Parameters are:  $L=800$, $\Delta t =0.5$, $T=60$, $N_s=16384$ with the renormalization factors obtained from system (\ref{Cons-mass-mom-a}) for $j$=1, 2. {Note the renormalized Duhamel iterations converge over such a large time-span, a significant improvement over when a single quantity is conserved.}}
    \label{mass-mom-cons-60}
\end{figure}

\begin{figure}
       \centering
       \includegraphics[width=60mm]{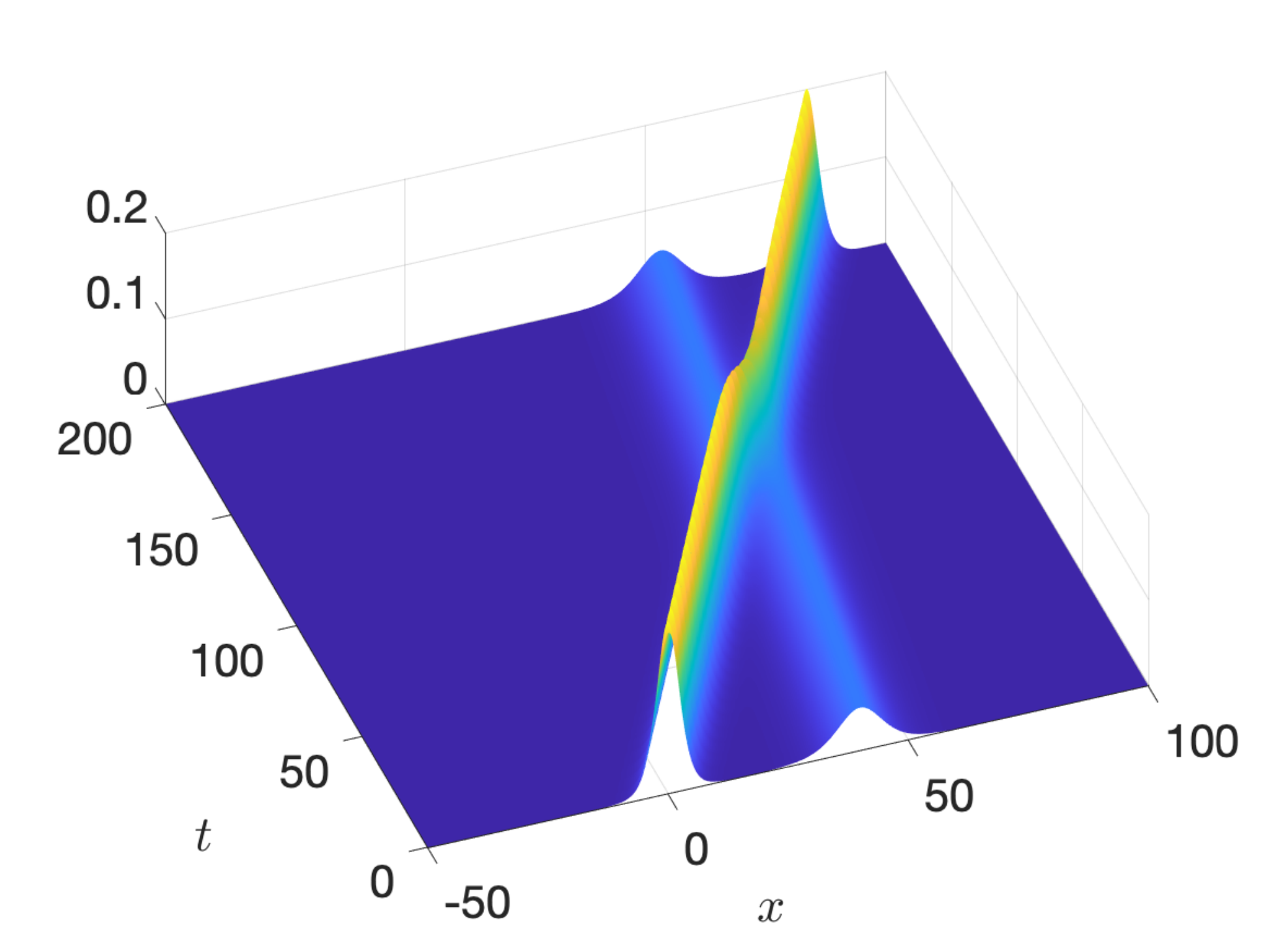}
       \caption{\small Interaction between two 1-soliton solution for the KdV equation with $\alpha =6, \epsilon = 1$ and initial condition
    $u(x,0)=(1/5){\rm sech}^2\left( x/\sqrt{10}\right)+(1/20){\rm sech}^2\left( (x-40)/2\sqrt{10}\right).$ Other parameters are $\Delta t = 0.5$, $L=800$, $N_S=16384$ and the time block size $T_1=20$. The renormalization factors are obtained by enforcing the conservation of mass and momentum simultaneously. The relative error in Hamiltonian $\approx7.86\times10^{-6}$ at $T=200$, while relative error in mass $\approx 1.9\times 10^{-16}$ and momentum $\approx 5.9\times 10^{-16}$ are kept near machine precision. }
       \label{2-soliton-figure-mass-mom}
   \end{figure}

%%%%%%%%%%%%%%%%%%%%%%%%%%%%%%%%%%%%%%%%%%%%%%%%%%%%%%%%%
Here the KdV solution is decomposed in the form: $u(x,t)= R_1(t) v_1(x,t)+R_2(t)v_2(x,t),$ where $R_1(t)$ and $R_2(t)$ are computed from the coupled system
%%%%%%%%%%%%%%%%%%%%%%%%%%%%%%%%%%%%%%%%%%%%%%
\begin{equation}
    \label{Cons-mass-mom-a}
    \mathcal{Q}_j\left[R_1(t)v_1(x,t) + R_2(t)v_2(x,t)\right]=\mathcal{Q}_j\left[ 2\beta^2{\rm sech}^2(\beta x) \right],\; j=1,2.
    \end{equation}
 %%%%%%%%%%%%%%%%%%%%%%%%%%%%%%%%%%%%%%%%%%%%%%
 {We choose $f_1(x)=(1/300){\rm sech}\left(x/\sqrt{600}\right)$, and  $f_2(x)=u_0(x)-f_1(x)$.}  The explicit expressions for $R_1(t)$ and $R_2(t)$ can be obtained from the coupled system
$$ R_1(t)= \frac{C_1-A_2(t)R_2(t)}{A_1(t)}, \quad \text{and} \quad \mu_1(t) R_2^2+\mu_2(t)R_2+ \mu_3(t)=0,$$
where 
\begin{equation*}
%\label{Renormalization-factors-mass-mom}
    %%%%%%%%%%%%%%%%%%%%%%%%%%%%%%%%%%%%%%%%%%%%%%%%%%
    \mu_1(t) = A_3 A_2^2+A_4A_1^2-2A_1 A_2 A_5, \quad %\nonumber
    %%%%%%%%%%%%%%%%%%%%%%%%%%%%%%%%%%%%%%%%%%%%%%%%%%
    \mu_2(t) = 2A_1 A_5 C_1-2A_2 A_3 C_1, \quad %\nonumber
    %%%%%%%%%%%%%%%%%%%%%%%%%%%%%%%%%%%%%%%%%%%%%%%%%%
    \mu_3(t) = A_3 C_1^2-C_2 A_1^2, %\nonumber
    \end{equation*} 
    and  
    
    \begin{eqnarray*}
    %%%%%%%%%%%%%%%%%%%%%%%%%%%%%%%%%%%%%%%%%%%%%%%%%%
    A_1(t) = \int_{\mathbb{R}} v_1(x,t) dx, \quad 
    A_2(t) &&= \int_{\mathbb{R}} v_2(x,t) dx, \quad A_3(t) = \int_{-\infty}^{\infty} v_1^2(x,t) dx, \\%\nonumber \\
    %%%%%%%%%%%%%%%%%%%%%%%%%%%%%%%%%%%%%%%%%%%%%%%%%%%%
   A_4(t) = \int_{\mathbb{R}} v_2^2(x,t) dx, \quad  A_5(t) &&= \int_{\mathbb{R}} v_1(x,t)v_2(x,t) dx.%\\%\nonumber \\
    %%%%%%%%%%%%%%%%%%%%%%%%%%%%%%%%%%%%%%%%%%%%%%%%%%
        %%%%%%%%%%%%%%%%%%%%%%%%%%%%%%%%%%%%%%%%%%%%%%%%%%
\end{eqnarray*}
%%%%%%%%%%%%%%%%%%%%%%%%%%%%%%%%%%%%%%%%%%%%%%%%%%%%%%%
Numerical tests indicate that the algorithm converges to the correct solution when using the root $R_2(t) = [-\mu_2(t)+\sqrt{\mu_2^2(t)-4\mu_1(t) \mu_3(t)}]/[2 \mu_1(t)]$, while the other one causes the TDSR algorithm to diverge. The correct root was found to always satisfy $R_2(0)v_2(x,0)=f_2(x)$ at any Duhamel iteration while the other consistently violated it.  Figure \ref{mass-mom-cons-60} shows results of TDSR simulations where mass and momentum are both conserved.
%%%%%%%%%%%%%%%%%%%%%%%%%%%%%%%%%%%%%%%%%%%%%%%%%%%%%%%%%%%%%%%%%%%%%%%%%%%%%%%%%%%%%%%%%%%%%%
%%%%%%%%%%%%%%%%%%%%%%%%%%%%%%%%%%%%%%%%%%%%%%%%%%%%%%%%%%%%%%%%%%%%%%
Another interesting numerical experiment (discussed below), is related to interaction (or collision) between two 1-soliton solutions to the KdV equation. The corresponding initial condition is
   \begin{equation}
   \label{initial-condition-1-soliton}
        u_0(x)=2\beta_1^2{\rm sech}^2\left(\beta_1 x\right)+2\beta_2^2{\rm sech}^2\left(\beta_2 (x-x_0)\right),
   \end{equation}
 with $\beta_1=\frac{1}{\sqrt{10}},\;\beta_2=\frac{1}{2\sqrt{10}}$ and an initial separation of $x_0=40$. We simulated this to end time $T=200$ using the idea of multi-blocking (see Fig. ~\ref{2-soliton-figure-mass-mom}). The solitons interact elastically, mainly emerging unscathed from the interaction, suffering only from a phase shift as expected. {As before, we prescribed $f_1(x)=(1/300){\rm sech}\left(x/\sqrt{600}\right)$, while $f_2(x)=u_0(x)-f_1(x)$}.
%%%%%%%%%%%%%%%%%%%%%%%%%%%%%%%%%%%%%%%%%%%%%%%%%%%%

\bigskip
%%%%%%%%%%%%%%%%%%%%%%%%%%%%%%%%%%%%%%%%%%%%%%%%%%%%%%%%%
\subsection{Conservation of mass and Hamiltonian}
\begin{figure}
    \centering
    \includegraphics[width=150mm]{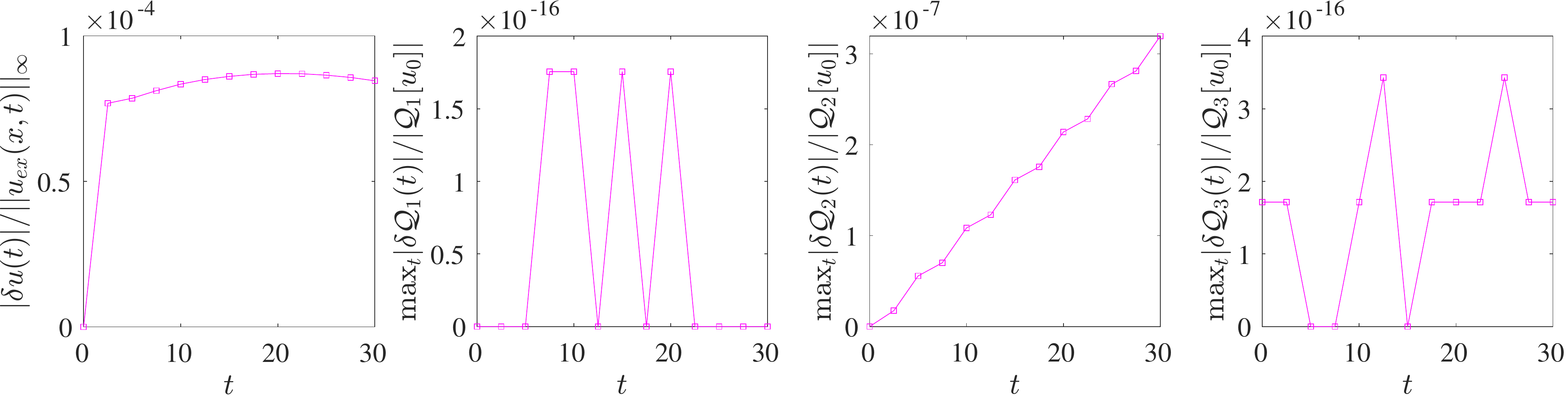}
    \caption{ \small (a) Time evolution of the relative error between the numerically obtained (TDSR) solution to the KdV equation with conservation of mass and Hamiltomian compared to its exact solution given in Eq.~(\ref{solution-accuracy}). Time evolution of the conserved quantities given in
Eq.~(\ref{cons-law-convergence}): mass (b), momentum (c) and Hamiltonian (d). Note that the relative errors in conservation of mass and Hamiltonian now stay close to machine precision. Parameters are:  $L=800$, $\Delta t =0.5$, $T=30$, $N_s=16384$ with the renormalization factors obtained from system (\ref{Cons-mass-mom-a}) for $j$=1, 3.}
    \label{mass-ham-cons-30}
\end{figure}

%%%%%%%%%%%%%%%%%%%%%%%%%%%%%%%%%%%%%%%%%%%%%%%%%%%%%%%%%%
In this case, the renormalization factors $R_1(t)$ and $R_2(t)$ are now computed from the coupled system of equations given in (\ref{Cons-mass-mom-a}) with $j$ taking the values $1$ and $3.$ {Like before, we pick the pseudo-initial conditions to be $f_1(x)=(1/300){\rm sech}\left(x/\sqrt{600}\right)$ and $f_2(x)=u_0(x)-f_1(x)$.} We solved system (\ref{Cons-mass-mom-a}) for $j=1,3$ using the Newton's method. 
Our findings are similar to those reported for the simultaneous conservation of mass and momentum, i.e., conservation of mass and Hamiltonian are achieved (see Fig. \ref{mass-ham-cons-30}). 
%%%%%%%%%%%%%%%%%%%%%%%%%%%%%%%%%%%%%%%%%%%%%%%%%%%%%%%%%%%%%%%%%%%

\bigskip
%%%%%%%%%%%%%%%%%%%%%%%%%%%%%%%%%%%%%%%%%%%%%%%%%%%%%%%
\subsection{Conservation of mass, momentum and Hamiltonian}
\begin{figure}
    \centering
    \includegraphics[width=150mm]{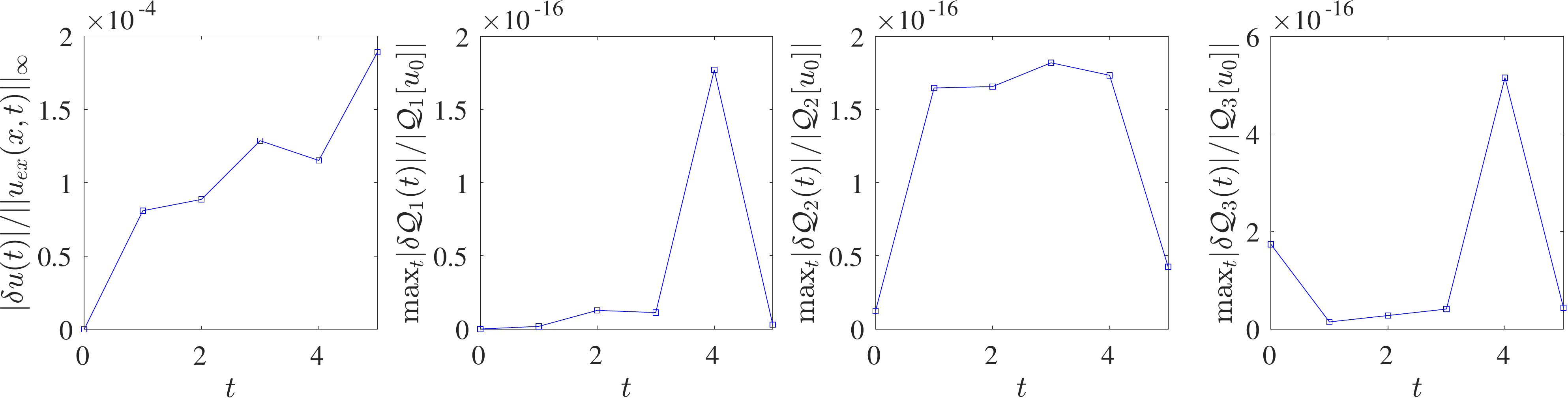}
    \caption{\small (a) Time evolution of the relative error between the numerically obtained solution, from TDSR with conservation of mass, momentum and Hamiltonian, for the KdV equation compared to its exact solution given in Eq.~(\ref{solution-accuracy}). Time evolution of the conserved quantities given in Eq.~(\ref{cons-law-convergence}): mass (b), momentum (c) and Hamiltonian (d). Parameters are: $\Delta t =0.5$, $T=5$, $N_s=2048$, $L=100$ with the renormalization factors obtained from Eqs.(\ref{Cons-mass-mom-ham-a}).}
    \label{Mass-ham-mom-convergence}
\end{figure}

%%%%%%%%%%%%%%%%%%%%%%%%%%%%%%%%%%%%%%%%%%%%%%%%%%%%%
Lastly, we seek three renormalization factors $R_1(t)$, $R_2(t)$ and $R_3(t)$ that satisfy $u(x,t)=\sum_{j=1}^{3}R_{j}(t)v_{j}(x,t)$ and obey the conservation laws 
\begin{equation}
    \label{Cons-mass-mom-ham-a}
    \mathcal{Q}_j\left[\sum_{\ell = 1}^3R_\ell (t)v_\ell (x,t)\right] = \mathcal{Q}_j\left[ 2\beta^2{\rm sech}^2(\beta x) \right],\;\;\; j=1,2,3.
    \end{equation}
 %%%%%%%%%%%%%%%%%%%%%%%%%%%%%%%%%%%%%%%%%%%%%%
 Unlike the previous cases, of mass-momentum and mass-Hamiltonian conservation, two of the current conserved quantities now are nonlinear functionals further limiting the choices for the pseudo-initial conditions as far as the convergence over large time intervals is concerned. {It turns out that convergence over moderately long time intervals ($T=5$) is achieved with pseudo-initial conditions in the Duhamel integral formulas Eqs.~(\ref{eqn-phi-1}) and (\ref{eqn-phi-j}) as 
 $f_3(x)=0.15\exp (-x^2), f_2(x)=0.05\exp (-x^2)$ and $f_1(x)=u_0(x)-f_2(x)-f_3(x)$.} 
%%%%%%%%%%%%%%%%%%%%%%%%%%
The renormalization factors are found from the coupled system of equations derived from enforcing the conservation of mass, momentum and conservation of Hamiltonian simultaneously Eq.~(\ref{Cons-mass-mom-a}) for $j=1,2,3$, using the Newton's root finding method. The results are depicted in Fig.~\ref{Mass-ham-mom-convergence} where the error in the momentum, Hamiltonian and mass are kept at the level of machine precision. 

\bigskip
\bigskip
%%%%%%%%%%%%%%%%%%%%%%%%%%%%%%%%%%%%%%%%%%%%%%%%%%%%%%%%%%%%%%%%%%%%%%%%%%%%%
\section{TDSR method: Multi-dimensional test case}
%%%%%%%%%%%%%%%%%%%%%%%%%%%%%%%%%%%%%%%%%%%%%%%%%%%%%%%%%%%%%%%%%%%%%%%%%%%%%%%%%%%%%%%%%
The nonlinear Schr\"odinger equation
\begin{eqnarray}
\label{2D-NLS}
iu_t+V({\bf x})u + \nabla^2 u + |u|^{2}u=0, 
\end{eqnarray}
plays an important role in modeling fundamental physics ranging from photonics, Bose-Einstein condensation to fluid mechanics \cite{ablowitz2011nonlinear}. Depending on the physics at hand, $\nabla^2 = \partial^2/\partial x^2 + \partial^2/\partial y^2$ denotes wave diffraction, $V$ is the refractive index, photonic lattice or an external potential 
and $|u|^2$ measuring the intensity or density of a complex valued wavefunction $u.$  As such, using the TDSR method to
simulate the NLS equation would seem natural. In the absence of any external potential ($V(\textbf{x})=0$), Eq.~ (\ref{2D-NLS}) admits a special class of solutions known as the Townes solitons. They are of the form $u(\textbf{x},t)=U_\lambda (\textbf{x})e^{i\lambda^2t}, \lambda\in\mathbb{R}$ with real valued function $U_\lambda$ satisfying the boundary value problem
%%%%%%%%%%%%%%%%%%%%%%%%%%%%%%%%%%
\begin{eqnarray}
\label{BVP-Townes-soliton}
 \nabla^2 U_\lambda + U^3_\lambda =  \lambda^2 U_\lambda.
\end{eqnarray}
%%%%%%%%%%%%%%%%%%%%%%%%%%%%%%%%%%
\begin{figure}
    \centering
    \includegraphics[width=75mm]{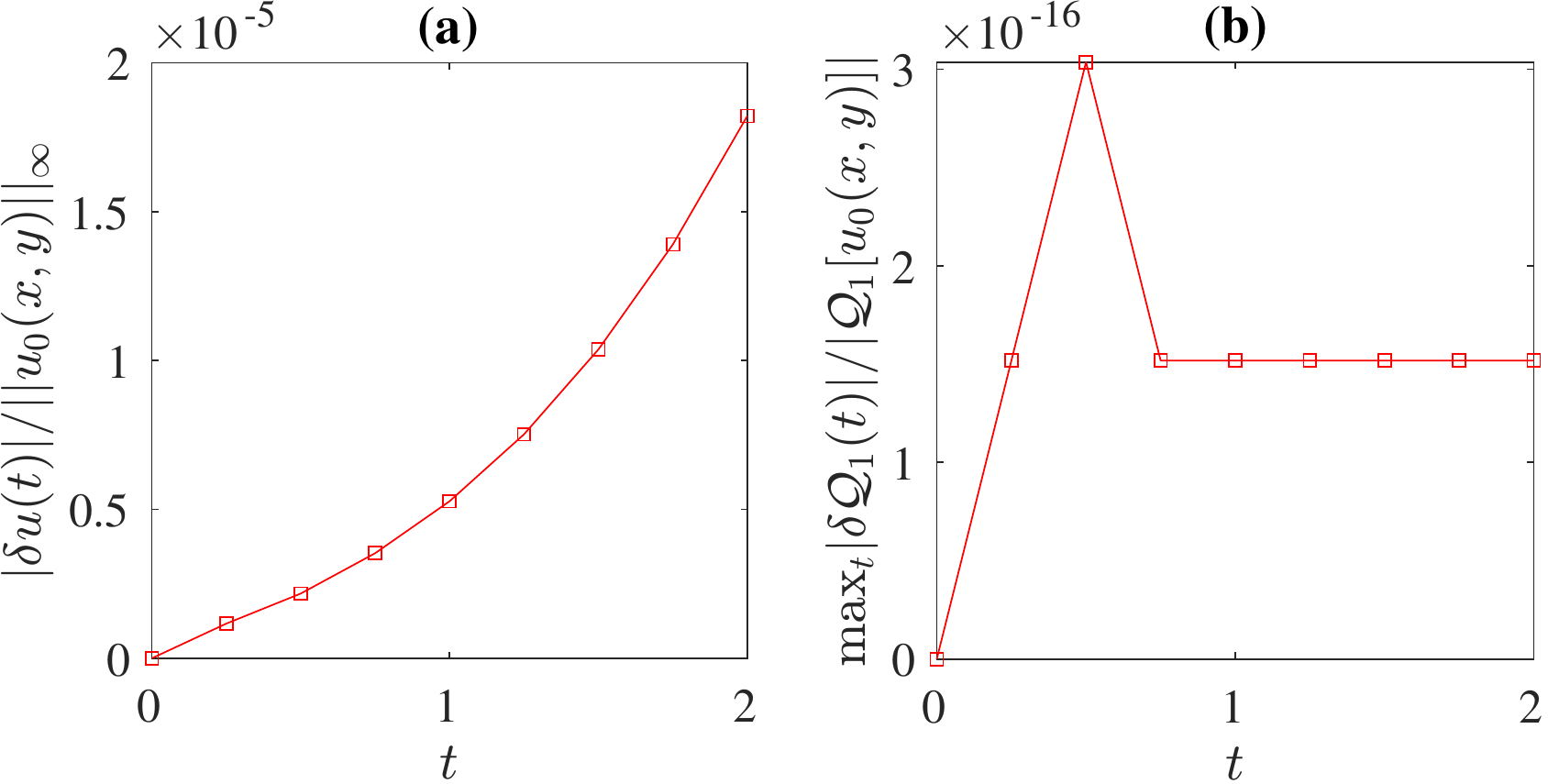}
    \caption{\small (a) {The error in the TDSR solution, as defined relative to the initial profile of the Townes soliton, (b) the relative error in the power as a function of time $t$. As one can see, it is kept close to machine precision over the entire time span [0,2]. Parameters are: $\Delta t=0.05$, a square spatial domain with $L=40$ , and a spatial discretization of $512\times512$  for the computations.}}
    \label{2D-NLS-error-plot}
\end{figure} 

The numerical computation of the Townes soliton is well documented in the literature and can be achieved with the use of various boundary value problem solvers (Quasi-Newton methods, spectral renormalization, etc. \cite{yang2009newton,Ablowitz:05,kevrekidis2007emergent,kivshar2003optical}).  While the implementation of the TDSR method on the NLS equation was first reported in \cite{cole2017time} it was exemplified in one spatial dimension for which the semi-group $\exp (it\partial^2/\partial x^2)$ lives in $(1+1)$ dimensions. Here, the $(1+2)D$ NLS equation, subject to sufficiently rapidly decaying boundary conditions is chosen as a prototypical example to demonstrate the applicability of the TDSR scheme in multiple spatial dimensions. 
The NLS Eq.~(\ref{2D-NLS}) admits three conservation laws: power $Q_1=|u|^2$, 
Hamiltonian $Q_3=-\frac{1}{4}|u|^4 +\frac{1}{2}|\nabla u|^2$ and momentum $Q_2=u\nabla u^{*}$.  We were able to reproduce the time evolution of the Townes soliton, to within a relative error of $\sim 1.8 \times 10^{-5}$, when the algorithm was implemented on time interval $[0, 2]$ using time step of size $\Delta t=0.05$ and number of spatial (square) grid points $N_{S}=512$. We rigorously conserve the power $Q_1$, and as can be seen in Fig.~\ref{2D-NLS-error-plot} (b), the relative error in power stays close to machine precision.

\bigskip
\bigskip 
%%%%%%%%%%%%%%%%%%%%%%%%%%%%%%%%%%%%%%%%%%%%%%%%%%%%%%%%%%%%
\section{Conclusions and future directions}
\label{sec-conclusions}
%%%%%%%%%%%%%%%%%%%%%%%%%%%%%%%%%%%%%%%%%%%%%%%%%%%%%%%%%%%%
In 2005 Ablowitz and Musslimani proposed the spectral renormalization method as a tool to numerically approximate solutions to nonlinear boundary value problems. Since then, it has been successfully used in many physical settings that include photonics \cite{yang2003fundamental}, Bose-Einstein condensation \cite{akkermans2008numerical}, Kohn-Sham density functional theory \cite{Kohn_Sham_Musslimani}, and water waves \cite{ablowitz2006new}.
%\af{These two paragraphs are not linked}
%It is usually the case that the method of choice for time-stepping evolution equations is fundamentally different than its time-independent counterpart. For example, computing steady state solution for the Navier-Stokes equations invokes 
%Newton-type solver as opposed to lower/higher order explicit-implicit Euler like numerical approximations of time dependent solutions.
%%%%%%%%%%%%%%%%%%%%%%%%%%%%%%%%%%%%%%%%%%%%%%%%%%%%%%%%%
In 2016, Cole and Musslimani proposed the time dependent spectral renormalization method to simulate evolution equations with periodic boundary conditions. This important idea brings two novel aspects: (i) it extends the original steady state spectral renormalization method to the time domain, thus offering a unifying approach by which time-independent as well as evolution equations are solved by the same numerical scheme, (ii) it allowed the inclusion of certain physics such as conservation and dissipation laws. 
In this paper we have significantly empowered the computational capabilities of the TDSR method that allows the (i) enforcement of several conservation laws or dissipation rate equations, (ii) flexibility to apply other non-periodic boundary conditions. We have successfully demonstrated these ideas on prototypical dynamical systems of physical significance. Examples include the Korteweg-de Vries equation and dynamics of fronts modeled by the
Allen-Cahn equation. We conclude this section by making a remark regarding possible application of the TDSR to weak wave turbulence.
Wave turbulence describes the chaotic interactions of dispersive wavetrains (analogs to eddies) when an external forcing term added to the underlying nonlinear evolution equations are mediated by dissipative forces (see \cite{newell2011wave,nazarenko2011wave,newell2001wave1} and references contained). Numerical investigations of these phenomena is a challenging task, requiring very long time runs for statistical equilibrium to be reached. Carefully designed numerical integrators are used (see \cite{majda1997one}) to integrate the potentially stiff, underlying nonlinear field equations over long time. In such scenarios, a close control over the conserved quantities of the dynamical system can prove vital to ensure long time accuracy of the solution. The application of the TDSR to this field thus seems natural and is kept for future work.
%%%%%%%%%%%%%%%%%%%%%%%%%%%%%%%%%%%%%%%%%%%%%%
\bigskip 
\bigskip

%%%%%%%%%%%%%%%%%%%%%%%%%%%%%%%%%%%%%%%%%%%%%%%%%%%%%%%%%%%
%\begin{acknowledgements}
%\end{acknowledgements}

%%%%%%%%%%%%%%%%%%%%%%%%%%%%%%%%%%%%%%%%%%%%%%%%%%%%

\end{document}